\documentclass[11pt]{amsart}
\usepackage[utf8]{inputenc}
\usepackage[T1]{fontenc}
\usepackage{lmodern}
\usepackage{amssymb,amsmath,amsthm,booktabs,enumitem,array}
\usepackage[margin=1.1in]{geometry}
\usepackage[colorlinks=true,linkcolor=blue,citecolor=blue,urlcolor=blue]{hyperref}
\hypersetup{
  pdftitle={First Neumann Eigenfunctions on Piecewise Smooth Jordan Domains with Convex Corners},
  pdfauthor={Santiago Pineda Montoya; Edisson Gallego Gonz\'alez; Jaiberth Porras Barrera; Johan H. R\'ua Mu\~noz}
}
\newtheorem{theorem}{Theorem}

\newtheorem{proposition}{Proposition}[section]
\newtheorem{lemma}[proposition]{Lemma}
\newtheorem{corollary}[proposition]{Corollary}
\theoremstyle{definition}
\newtheorem{definition}[proposition]{Definition}
\theoremstyle{remark}
\newtheorem{remark}[proposition]{Remark}
\numberwithin{equation}{section}
\newcommand{\R}{\mathbb{R}}
\newcommand{\D}{\mathbb{D}}
\newcommand{\Sone}{\mathbb{S}^1}
\newcommand{\cG}{\mathcal{G}}
\newcommand{\cP}{\mathcal{P}}
\newcommand{\cA}{\mathcal{A}}
\newcommand{\cH}{\mathcal{H}}
\newcommand{\Gq}{G_q}
\newcommand{\Kq}{K_q}
\newcommand{\Nq}{N_q}
\newcommand{\dv}{\operatorname{div}}
\newcommand{\curl}{\operatorname{curl}}
\newcommand{\ip}[2]{\langle #1,#2\rangle}
\newcommand{\doi}[1]{\href{https://doi.org/#1}{doi:\nolinkurl{#1}}}
\title[Piecewise smooth domains with convex corners]{First Neumann eigenfunctions on piecewise smooth\\ Jordan domains with convex corners}
\author[S. Pineda Montoya]{Santiago Pineda Montoya}
\address{Universidad Nacional de Colombia, Medell\'in, Colombia}
\email{sapinedamo@unal.edu.co}
\thanks{Corresponding author: Santiago Pineda Montoya (\texttt{sapinedamo@unal.edu.co}).}

\author[E. Gallego Gonz\'alez]{Edisson Gallego Gonz\'alez}
\address{Universidad EAFIT, Medell\'in, Colombia}
\email{egalleg8@eafit.edu.co}

\author[J. Porras Barrera]{Jaiberth Porras Barrera}
\address{Universidad Nacional de Colombia, Leticia, Amazonas, Colombia}
\email{jporras@unal.edu.co}

\author[J. H. R\'ua Mu\~noz]{Johan H. R\'ua Mu\~noz}
\address{Universidad de Antioquia, Medell\'in, Colombia}
\email{heriberto.rua@udea.edu.co}

\date{19 November 2026}
\begin{document}
\begin{abstract}
We prove that every nonzero eigenfunction of the first positive Neumann
eigenvalue on a bounded planar Jordan domain whose boundary consists of finitely
many smooth arcs meeting at interior angles strictly between $0$ and $\pi$ has
nonvanishing gradient in the interior; hence its extrema are attained only on
the boundary. The eigenvalue may be multiple and the domain need not be convex;
every convex polygon is included. The method adapts a recent Green-kernel
approach for smooth domains.
\end{abstract}
\maketitle

\section{Introduction and statement}

Rauch's hot spots conjecture \cite{Rau75} asserts that an eigenfunction of the
first positive Neumann eigenvalue of a planar domain attains its maximum and
minimum only on the boundary; see \cite{BB99} for precise formulations. It fails
for some domains with holes \cite{BW99}, and in high dimensions even convex
domains can have interior hot spots \cite{dDP24}. In the plane it is known, among
other cases, for domains with two axes of symmetry \cite{JN00}, for lip
domains \cite{AB04}, for triangles \cite{JM20,JM22,Siu15,CGY26}, for isosceles
trapezoids and kites \cite{DGJYYZ26}, for L-shaped domains and several
classes of domains tiled by them \cite{Hat25}, and for the union of a regular
polygon with its reflection across an edge \cite{Hat25b}; the last two families
include domains with reentrant corners.

A recent preprint \cite{OAI} asserts the strict form, absence of interior
critical points, for every smooth simply connected planar domain. Assuming that
assertion, Burdzy \cite{Bur26} and Ba\~nuelos \cite{Ban26} derive, by
approximation, the non-strict conclusion on nonsmooth simply connected domains
whose first positive eigenvalue is simple, and Ba\~nuelos asks whether it holds
on every simply connected polygon, the eigenvalue being possibly multiple
\cite[Question~1]{Ban26}. Independently of \cite{OAI}, Deng, Jiang and Yang
\cite{DJY26} prove that on every bounded convex planar domain the maximum of
$|u|$ is attained only on the boundary. This note treats domains with corners. It
does not assume the assertion of \cite{OAI}, and for convex polygons it answers
the question of \cite{Ban26} in the strict form.

All functions and Hilbert spaces in this note are real unless complex
coordinates are explicitly used for conformal maps. We use the standard Sobolev
spaces $W^{k,p}$, write $H^k=W^{k,2}$, and let $H^1_0$ denote the closure of
compactly supported smooth functions in $H^1$. The boundary space $H^{1/2}$
is the Sobolev trace space; $\operatorname{Tr}$ denotes the trace operator.
The notation $C^{k,\gamma}$ means that the derivatives of order $k$ are
H\"older continuous with exponent $\gamma$. These conventions apply componentwise
to vector fields. The symbol $B(x,r)$ denotes an open Euclidean disk,
$\delta_y$ the unit point mass at $y$.
The constants $C$ in estimates may change from line to line; $f\asymp g$
means two-sided comparison by positive constants in the indicated neighborhood.

A bounded Jordan domain means the bounded component enclosed by a simple closed curve. We use the following precise meaning of a piecewise smooth domain.

\begin{definition}\label{def:domain}
An admissible domain is a bounded planar Jordan domain $\Omega$ with a finite
set of corners $\{v_1,\ldots,v_m\}$, where $m\ge0$. Each component of the boundary with the corners
removed is a regular $C^\infty$ arc, smooth with nonzero one-sided tangent up to
its endpoints. At $v_j$ the interior angle is $\alpha_j\pi$, where
$0<\alpha_j<1$. The case of a $C^\infty$ Jordan boundary with no corners is also
allowed. Global convexity is not required. The open smooth part of the boundary
is denoted by $\Gamma$.
\end{definition}

In particular $\Omega$ is simply connected and Lipschitz; $|\Omega|$ denotes its area. Its outward unit normal
on $\Gamma$ is $\nu$, boundary arclength is $d\sigma$, and $\tau=(-\nu_2,\nu_1)$ gives the counterclockwise
orientation. A boundary vector is tangent when its scalar product with $\nu$ is zero; for Sobolev fields this condition is imposed on the trace almost everywhere on $\Gamma$. Signed curvature is fixed by $\partial_\tau\tau=-\kappa\nu$;
it is positive on a circle bounding a disk. It is bounded on the finite arcs,
but may have either sign. Put
\[
\mu=\inf\Big\{\int_\Omega|\nabla v|^2\Big/\int_\Omega v^2:
0\ne v\in H^1(\Omega),\ \int_\Omega v=0\Big\},
\]
and define its eigenspace by the weak formulation
\[
V=\left\{v\in H^1(\Omega):
\int_\Omega\nabla v\cdot\nabla\psi=\mu\int_\Omega v\psi
\text{ for every }\psi\in H^1(\Omega)\right\}.
\]
Write $\lambda_D$ for the first Dirichlet eigenvalue, defined by the same
Rayleigh quotient over $0\ne v\in H^1_0(\Omega)$, without the mean-zero constraint.
Poincar\'e's inequalities and compact Sobolev embedding on this connected
bounded Lipschitz domain give $0<\mu,\lambda_D<\infty$ and a finite-dimensional
eigenspace $V$.
The operator $-\Delta_D$ below is the Dirichlet realization of the Laplacian
on $L^2(\Omega)$. Unsubscripted $L^2$ norms in physical-domain arguments are
taken over $\Omega$. Euclidean vector products are denoted by a dot, and
$\nabla X:\nabla Y=\sum_{i,j}(\partial_iX_j)(\partial_iY_j)$ denotes the scalar product of the two Jacobian matrices.

\subsection*{Disk kernels} Let $\D=\{z\in\mathbb C:|z|<1\}$ be the unit disk,
identified with a subset of $\mathbb R^2$, and let $\Sone=\partial\D$ with
arclength $ds$. Its outward normal is $\nu_\D(s)=s$. Set
\[
G(x,y)=\frac1{4\pi}\log\Big(1+\frac{(1-|x|^2)(1-|y|^2)}{|x-y|^2}\Big),\qquad
P_s(x)=\frac{1-|x|^2}{2\pi|s-x|^2},\qquad N_0(s,t)=\frac1{\pi|s-t|^2}
\]
for $x,y\in\D$, $x\ne y$, and $s,t\in\Sone$, $s\ne t$.
These are respectively the Dirichlet Green function, the Poisson kernel and
the kernel of the Dirichlet integral. In nonnegative kernel integrals use
$G(x,x)=+\infty$; its diagonal is a null set for the absolutely continuous
measures considered here. Diagonal values of boundary kernels in double
arclength integrals are immaterial and may be set to zero. For a measure $d\rho=a_\rho\,dx$ on $\D$ with $a_\rho\ge 0$ measurable, and a nonnegative measurable
$f$, write $(T_\rho f)(x)=\int_\D G(x,z)f(z)\,d\rho(z)\in[0,\infty]$, and define, for $p\in\D$ and
$s\ne t$ in $\Sone$ (with $T_\rho^0=I$, the identity operator),
\begin{equation}\label{eq:kernels}
K^\rho(p,s)=\sum_{n\ge 0}(T_\rho^nP_s)(p),\qquad
R^\rho(s,t)=\sum_{n\ge 0}\int_\D P_s\,T_\rho^nP_t\,d\rho,\qquad N^\rho=N_0+R^\rho .
\end{equation}
Here $T_\rho^nP_s$ denotes an iterated nonnegative kernel integral, even
when $P_s$ does not belong to $L^2(\rho)$. Products involving an infinite
value and a zero nonnegative factor are assigned value zero.
These are sums of nonnegative integrals: $K^\rho,N^\rho$ have values in $(0,\infty]$ and $R^\rho$ in $[0,\infty]$; for bounded $a_\rho$
they are the kernels (4.6)--(4.8) of \cite{OAI}. A symmetric kernel $D$ with zero diagonal on a set
$S$ is \emph{conditionally negative semidefinite} if $\sum_{i,j}\omega_i\omega_jD(s_i,s_j)\le 0$ for all
finite families $s_i\in S$ and real $\omega_i$ with $\sum_i\omega_i=0$.

The reciprocal quotient in these kernels is the central analytic object.
The kernel assertion needed here concerns bounded nonnegative densities with
compact support in $\D$ and a strict Dirichlet form bound. We prove that
assertion in the kernel section before applying it to the conformal truncations;
the main statement has no additional reciprocal-kernel hypothesis.

\begin{theorem}\label{thm:A}
Let $\Omega$ be an admissible domain in Definition~\ref{def:domain}. Then for every $0\ne u\in V$,
\[
\nabla u(x)\ne 0\qquad(x\in\Omega),\qquad\text{and hence}\qquad
\min_{\partial\Omega}u<u(x)<\max_{\partial\Omega}u\qquad(x\in\Omega).
\]
\end{theorem}

A prevertex is a boundary point of the disk mapped to a corner by a conformal
bijection $\Phi:\D\to\Omega$.

\begin{remark}\label{rem:smooth}
If $m=0$, that is, for a $C^\infty$ Jordan boundary, Theorem~\ref{thm:A} is
\cite[Theorem~1.1]{OAI}, and the proof below follows that source except for
the finite preservation step; no priority is claimed in that case. The extension addressed here concerns corners: the conformal potential is singular at the
prevertices, the boundary gradient extends continuously and vanishes at the corners,
and the Schoenberg map is injective only off the prevertices. Lipschitz
regularity of the gradient is not assumed. The signed curvature term is
retained until it cancels.
\end{remark}

\subsection*{Strategy of the proof}
Fix $0\ne u\in V$ and suppose that $\nabla u(\xi_0)=0$ at an interior point
$\xi_0=\Phi(p)$, where $\Phi:\D\to\Omega$ is a conformal bijection. The argument
has four steps.

\emph{1. A nonnegative boundary energy with known null vectors.} For a tangent
vector field $X$ on $\Omega$, Rohleder's principle
(Proposition~\ref{prop:rohleder}) gives $Q(X)\ge\mu\|X\|_2^2$, with equality
exactly for gradients of eigenfunctions. Pulling back by $\Phi$ and extending
boundary data by solutions of $-\Delta W=qW$ with $q=\mu|\Phi'|^2$ turns
$Q-\mu\|\cdot\|_2^2$ into a nonnegative quadratic form $E$ on tangent boundary
data on $\Sone$ (Section~\ref{sec:energy}). The boundary gradient
$g=(\nabla u)\circ\Phi$ satisfies $E(g)=0$.

\emph{2. A multiplier identity.} For every Lipschitz scalar function $b$ on
$\Sone$,
\[
E(bg)=\tfrac12\iint_{\Sone\times\Sone}N(s,t)\,(b(s)-b(t))^2\,g(s)\cdot g(t)\,ds\,dt
\]
(Proposition~\ref{prop:mult}); the signed curvature terms cancel.

\emph{3. A squared-distance kernel.} For each $p\in\D$ the kernel
$D_p(s,t)=K(p,s)K(p,t)/N(s,t)$ is conditionally negative semidefinite
(Propositions~\ref{prop:compactkernel} and~\ref{prop:cnd}). Hence
$D_p(s,t)=\|\mathcal B(s)-\mathcal B(t)\|^2$ for a Lipschitz map $\mathcal B$
into a Hilbert space (Lemma~\ref{lem:embed}). Applying step~2 to the coordinates
$b_j$ of $\mathcal B$ and summing gives
\[
\sum_jE(b_jg)=\tfrac12\Big|\int_{\Sone}K(p,s)\,g(s)\,ds\Big|^2
=\tfrac12|\nabla u(\xi_0)|^2=0 .
\]

\emph{4. Contradiction.} Every $E(b_jg)$ vanishes, so every $b_jg$ is the
boundary gradient of an eigenfunction. Since $\mathcal B$ is injective off the
prevertices, this produces either three linearly independent eigenfunctions,
contradicting $\dim V\le2$ (Proposition~\ref{prop:dim}), or a nonzero
eigenfunction with vanishing Cauchy data on a boundary arc, contradicting unique
continuation.

\emph{Where the corners enter.} For a smooth boundary these are the steps of
\cite{OAI}. At a corner of angle $\alpha_j\pi$ the potential $q$ blows up like
$|z-s_j|^{2\alpha_j-2}$ at the prevertex. The kernel theorem is therefore proved
first for bounded compactly supported potentials and then passed to $q$ by
monotone truncation; this requires $K(p,\cdot)$ bounded and $N$ finite off the
prevertices (Corollary~\ref{cor:greenC1} and Lemma~\ref{lem:bdry}). The boundary
gradient $g$ is continuous and vanishes at the prevertices, but it need not
be Lipschitz (Lemma~\ref{lem:gradient}); this is why the energy is defined for
bounded $H^{1/2}$ data. The curvature identity must hold without vertex terms
(Lemma~\ref{lem:det}), and the topological step must allow $\mathcal B$ to fail
to be injective at finitely many points (Lemma~\ref{lem:top}). The condition
$\alpha_j<1$ is used exactly in these estimates.

\emph{Organization.} Section~\ref{sec:prelim} proves the corner estimates and
$\dim V\le2$; Section~\ref{sec:rohleder} the variational principle and the
curvature identity; Section~\ref{sec:kernels} the kernel theorem and the singular
limit; Section~\ref{sec:energy} the energy formula and the multiplier identity;
Section~\ref{sec:proof} the conclusion.

\subsection*{Relation to previous work}
The strategy, multiplicity argument, vector variational principle,
boundary kernels and Hilbert-space contradiction are adapted from \cite{OAI}.
The reciprocal-kernel argument is supplied for
bounded compactly supported densities. Its initial regularization and limiting
framework follow \cite[Section~5]{OAI}, with proofs included. The finite
preservation step is proved by a matrix flow and a cone-invariance argument.
This is an alternative to the rank-one update lemma \cite[Lemma~5.4]{OAI},
which preserves the reciprocal condition for a single column update and,
when iterated, gives the finite resolvent assertion; it is
the inclusion of a proof, not the choice of the flow, that removes
\cite[Theorem~5.1]{OAI} as an assumption. The regularization and limiting framework are shared with that source.
Their proofs are included here and must be checked as part of the present
argument; no assertion of Section~5 is used without proof.
Attribution of the statement and of the adapted constructions is retained. Monotone truncation also occurs in
that source; its use here needs the corner bounds proved below.

Straight convex polygons are special cases. Compared with that case, three
changes are essential: sector separation of variables is replaced by estimates
on the disk, tangency-preserving smoothing is carried out in curved coordinates,
and a signed curvature integral is retained until it cancels.
The variational principle and the treatment of curvature on such domains
are already present in \cite[Theorem~1.1 and Section~3]{Roh21}; we include
the proofs needed here and make no novelty claim for that framework. Boundary uniqueness
on a smooth arc uses extension by zero, not reflection in a line. These changes
are proved rather than inferred by analogy. The finite-exceptional-set topology
argument is still required.

The classical inputs are the spectral inequality of Filonov, homogeneous
Dirichlet and Neumann $H^2$ regularity for curvilinear polygons without reentrant
corners, boundary regularity of conformal maps, smooth-boundary elliptic
regularity, trace and Sobolev theorems, and Schoenberg's theorem. We spell out
their uses. In particular, the regularity input is local in the corner angles;
it does not require the whole domain to be convex.

For context, Rohleder \cite[Theorem~1.1 and Remark~2.2]{Roh26} gives unconditional
restrictions on interior critical points in convex domains, including absence
under the additional condition
$\mu\operatorname{diam}(\Omega)^2\le j_{1,1}^2$, where $j_{1,1}$ is the first positive zero of the Bessel function $J_1$ of the first kind and $\operatorname{diam}(\Omega)=\sup_{x,y\in\Omega}|x-y|$.
Ba\~nuelos and Pang \cite[Theorem~1.1]{BP06} prove finiteness of the connected
components of positive level sets under smoothness and strictly positive
curvature assumptions. These are comparisons of scope, not inputs to the present
proof.

Reentrant corners, cusps, infinitely many corners and rough arcs are outside
Definition~\ref{def:domain}. Nor is an approximation argument used to infer strict
nonvanishing from a non-strict limiting conclusion.

\section{Corner geometry and analytical estimates}\label{sec:prelim}

We first separate the classical regularity facts from two estimates that will
replace the explicit Bessel expansions available for straight corners.

\begin{lemma}[classical inputs]\label{lem:std}
\begin{enumerate}[label=\textup{(\alph*)},leftmargin=2em]
\item\label{it:gap} $\mu<\lambda_D$.
\item\label{it:H2} For $f\in L^2(\Omega)$, the weak Dirichlet solution of
$-\Delta\varphi_D=f$, $\varphi_D\in H^1_0(\Omega)$, and the weak homogeneous Neumann solution of
$-\Delta\varphi_N+\varphi_N=f$ belong to $H^2(\Omega)$. In particular $V\subset H^2(\Omega)$,
and its elements are bounded and continuous on $\overline\Omega$.
\item\label{it:conf} A conformal bijection $\Phi:\D\to\Omega$ extends to a
homeomorphism of the closures. Put $s_j=\Phi^{-1}(v_j)$ and
$A=\Sone\setminus\{s_1,\ldots,s_m\}$. On every compact subarc of $A$, $\Phi$ is
smooth up to the boundary and $\Phi'\ne0$. Near a prevertex,
\begin{equation}\label{eq:SC}
|\Phi'(z)|\asymp |z-s_j|^{\alpha_j-1},\qquad
|\Phi(z)-v_j|\asymp |z-s_j|^{\alpha_j}.
\end{equation}
The comparisons hold for approach through $\overline\D\setminus\{s_j\}$.
\item\label{it:green} The Dirichlet operator $-\Delta-\mu$ on $\Omega$ has a
positive symmetric Green function $\cG$. With $G_0^\Omega$ the Laplace Green function,
\begin{equation}\label{eq:resOmega}
\cG(\xi,\eta)=G_0^\Omega(\xi,\eta)
+\mu\int_\Omega G_0^\Omega(\xi,\zeta)\cG(\zeta,\eta)\,d\zeta.
\end{equation}
Its $L^2(\Omega)$ sections are bounded uniformly in the pole. Away from the
diagonal it is jointly smooth up to the open smooth arcs in either variable,
and vanishes there in each variable.
\end{enumerate}
\end{lemma}

\begin{proof}[References and explanations]
Part~\ref{it:gap} is Filonov's inequality \cite{Fil05} for bounded Lipschitz domains.
For \ref{it:H2}, use the homogeneous Dirichlet and Neumann regularity theory
of \cite{Gri85} for curvilinear polygons whose individual arcs are of class
$C^{1,1}$ up to their endpoints, as in Definition~\ref{def:domain}.
The whole boundary is not asserted to be $C^{1,1}$ across a corner.
A compact statement for such polygons, taken from \cite{Gri85}, is collected in
\cite[Appendix~A.1, Theorem~A.2]{Cri18}. For either pure boundary condition the
first nonconstant corner exponent is $\pi/(\alpha_j\pi)=1/\alpha_j>1$.
Thus no singular term below the $H^2$ threshold occurs. Smooth local boundary
charts handle the rest of $\partial\Omega$. This is an angular regularity
criterion, not an appeal to global convexity. The Sobolev embedding
$H^2(\Omega)\subset C^{0,\gamma}(\overline\Omega)$ for some $\gamma>0$ gives the
last assertion.

Part~\ref{it:conf} uses Carath\'eodory's theorem and local Kellogg--Warschawski
regularity; see \cite[Chapters~2--3]{Pom92}. The corner derivative estimates are
the usual local power-map consequence of this regularity; see also the proof of
\cite[Lemma~4.3, equations~(8)--(10)]{PP}. Specifically, the power
$(w-v_j)^{1/\alpha_j}$, with a local branch on the interior side of the corner,
opens its angle to $\pi$. The two smooth arcs become a $C^{1,\delta}$ boundary
at the image point, for some $\delta>0$. The normalized local conformal map then
has a nonzero continuous derivative there. Factoring out the power gives
\eqref{eq:SC}. No holomorphic extension through a merely smooth arc is asserted.

For \ref{it:green}, the Dirichlet resolvent exists by \ref{it:gap}, and
$(-\Delta_D-\mu)^{-1}=\sum_{k\ge0}\mu^k(-\Delta_D)^{-k-1}$ is positivity
preserving. Set
$\cG(\cdot,\eta)=G_0^\Omega(\cdot,\eta)
+\mu(-\Delta_D-\mu)^{-1}G_0^\Omega(\cdot,\eta)$.
For any fixed $R_*>\operatorname{diam}(\Omega)$, domain monotonicity and
the Green function of the disk $B(\eta,R_*)$ give
$0\le G_0^\Omega(\xi,\eta)\le (2\pi)^{-1}\log(R_*/|\xi-\eta|)$.
Its square has a uniformly bounded integral over that disk. Together with the
bounded resolvent $L^2\to H^2\to L^\infty$, this proves uniform $L^2$
sections and $0\le\cG-G_0^\Omega\le C$. The resolvent identity gives \eqref{eq:resOmega}.
Local smooth-boundary elliptic regularity, applied in each variable on disjoint
neighborhoods, gives the asserted joint regularity off the diagonal.
\end{proof}

Fix once and for all a map $\Phi$ as in Lemma~\ref{lem:std}\ref{it:conf}, together with its prevertices $s_j$ and the regular boundary set $A$. Write $a=\mu|\Phi'|^2$. By \eqref{eq:SC}, for some $1<p_0<\infty$,
\begin{equation}\label{eq:apotential}
a\in L^{p_0}(\D),\qquad
0\le a(z)\le C\left(1+\sum_{j=1}^m |z-s_j|^{2\alpha_j-2}\right).
\end{equation}
When $m\ge1$, one may choose $1<p_0<\min_j(1-\alpha_j)^{-1}$. If there are no corners,
$a$ is bounded and all the following boundary estimates are immediate from
smooth-boundary regularity.

We record the Sobolev changes of variables used below. Composition with
$\Phi$ is an isometry for the Dirichlet seminorm on compactly supported smooth
functions. Completing these spaces in that seminorm, which is equivalent to
the $H^1_0$ norm by Poincar\'e's inequality on each bounded domain, gives
an isomorphism between $H^1_0(\Omega)$ and $H^1_0(\D)$. The limits agree
with ordinary composition on compact subsets, where $\Phi$ is a smooth
diffeomorphism. For a general $\psi\in H^1(\D)$, put
$p_0'=p_0/(p_0-1)$. H\"older's inequality gives
\[
\int_\Omega|\psi\circ\Phi^{-1}|^2
=\mu^{-1}\int_\D a|\psi|^2
\le\mu^{-1}\|a\|_{L^{p_0}}\|\psi\|_{L^{2p_0'}}^2<\infty.
\]
The local chain rule and the invariant Dirichlet integral thus place
$\psi\circ\Phi^{-1}$ in $H^1(\Omega)$. Conversely a bounded physical
$H^1$ function pulls back to $H^1(\D)$, by the same energy identity and
boundedness on the disk. These facts do not use any boundary kernel.

\begin{lemma}[vanishing of the physical gradient]\label{lem:gradient}
For $u\in V$, $\nabla u\in H^1(\Omega;\R^2)\cap C(\overline\Omega;\R^2)$ and
$\nabla u(v_j)=0$. More precisely, with $r=|z-s_j|$ small,
\begin{equation}\label{eq:gradcorner}
|\nabla u(\Phi(z))|\le C
\begin{cases}
r^{\alpha_j},&0<\alpha_j<\tfrac12,\\
r^{1/2}(1+|\log r|),&\alpha_j=\tfrac12,\\
r^{1-\alpha_j},&\tfrac12<\alpha_j<1.
\end{cases}
\end{equation}
On $\Gamma$, the gradient is smooth and tangent.
\end{lemma}

\begin{proof}
By Lemma~\ref{lem:std}, $u\in H^2$ is bounded and smooth up to $\Gamma$.
Put $U=u\circ\Phi$. Conformal invariance gives $U\in H^1(\D)\cap L^\infty(\D)$ and
\begin{equation}\label{eq:Neudisk}
\int_\D\nabla U\cdot\nabla\psi=\int_\D aU\psi
\quad(\psi\in H^1(\D)),\qquad \int_\D aU=0.
\end{equation}
To justify the full test space, $\psi\circ\Phi^{-1}$ has finite Dirichlet
integral and belongs to $L^2(\Omega)$: by \eqref{eq:apotential} and the planar
Sobolev embedding, $\int_\D a\psi^2<\infty$. It is consequently an admissible
test in the physical Neumann problem.

The normalized Neumann Green function of the disk is
\[
G_N(z,y)=-\frac{\log|z-y|+\log|1-z\overline y|}{2\pi}
+\frac{|z|^2+|y|^2}{4\pi}-\frac3{8\pi}.
\]
Direct differentiation and integration verify
$-\Delta_zG_N=\delta_y-1/\pi$, $\partial_{\nu_\D}G_N=0$ and
$\int_\D G_N(z,y)\,dz=0$. Furthermore,
\[
|1-z\overline y|^2-|z-y|^2=(1-|z|^2)(1-|y|^2)\ge0,
\qquad |\nabla_zG_N(z,y)|\le C(1+|z-y|^{-1}).
\]
The representation of the mean-zero part of $U$ with source $aU\in L^{p_0}$
therefore gives, at interior points,
\begin{equation}\label{eq:gradU}
|\nabla U(z)|\le C\int_\D(1+|z-y|^{-1})a(y)\,dy.
\end{equation}
One can obtain the representation first for smooth mean-zero sources and then
pass to $L^{p_0}$ using the Neumann solution operator on the smooth disk. For
the gradient formula at a fixed interior point, the source is smooth locally;
the remaining integral is differentiated away from its singularity.

Fix $0<r_0<1$. For $\beta>0$, $s\in\Sone$ and $r=|z-s|<r_0/4$, a direct splitting of the integral gives
\begin{equation}\label{eq:convolution}
\int_{\D\cap B(s,r_0)}\frac{|y-s|^{\beta-2}}{|z-y|}\,dy
\le C\begin{cases}
r^{\beta-1},&0<\beta<1,\\
1+|\log r|,&\beta=1,\\
1,&\beta>1.
\end{cases}
\end{equation}
Here is the estimate explicitly. On $|y-s|<r/2$, the denominator is at least
$r/2$, giving $Cr^{\beta-1}$. On $r/2\le|y-s|\le2r$, the weight is at most
$Cr^{\beta-2}$ and the integral of $|z-y|^{-1}$ over that set is at most $Cr$.
On $2r<|y-s|<r_0$, the denominator is at least $|y-s|/2$, giving
$C\int_{2r}^{r_0}t^{\beta-2}\,dt$. These are the three cases in
\eqref{eq:convolution}; the other prevertices contribute a bounded amount.

Apply this with $\beta=2\alpha_j$ in \eqref{eq:gradU}, and use
$|\nabla u(\Phi(z))|=|\nabla U(z)|/|\Phi'(z)|$ and \eqref{eq:SC}.
This gives \eqref{eq:gradcorner}, whose right side tends to zero in every case.
Together with smoothness up to each open arc this gives the continuous extension
of the gradient, while its $H^1$ membership follows from $u\in H^2$.
\end{proof}

\begin{lemma}[a Dirichlet regularity iteration]\label{lem:bootstrap}
Suppose $a$ obeys \eqref{eq:apotential}. If $v\in H^1_0(\D)\cap L^\infty(\D)$
satisfies $-\Delta v=av+f$ with $f\in L^\infty(\D)$, then
$v\in W^{2,p}(\D)$ for some $p>2$. In particular
$v\in C^{1,\gamma}(\overline\D)$ for some $\gamma>0$.
\end{lemma}

\begin{proof}
The case $m=0$ follows from smooth-disk regularity with bounded $a$. Otherwise let $\alpha_* =\min_{1\le j\le m}\alpha_j>0$.
Initially $av+f\in L^{p_0}$ for any $1<p_0<\min_j(1-\alpha_j)^{-1}$. Varying $p_0$, Dirichlet $W^{2,p_0}$
regularity on the disk and Sobolev embedding give
$v\in C^{0,\gamma_0}(\overline\D)$ for every
$0<\gamma_0<\min(1,2\alpha_*)$; only such subunit exponents will be used.
Since the boundary trace is zero, $|v(z)|\le C|z-s_j|^{\gamma_0}$ near every $s_j$.
Consequently, if such an exponent $\gamma_0\in[0,1)$ is available, then
\[
|av+f|\le C\left(1+\sum_j |z-s_j|^{2\alpha_j-2+\gamma_0}\right).
\]
This belongs to $L^p$ whenever $p>1$ and
$p(2-2\alpha_*-\gamma_0)<2$ if the expression in parentheses is positive;
if it is nonpositive there is no finite upper restriction on $p$.
Another $W^{2,p}$ estimate therefore yields every exponent
\[
0<\gamma_1<\min(1,\gamma_0+2\alpha_*).
\]
Starting with boundedness, finitely many iterations give an exponent
$0<\gamma_0<1$ with $\gamma_0>1-2\alpha_*$. For example one can increase the
exponent in steps smaller than $\alpha_*$ until crossing that threshold,
choosing a last step below $1$. The displayed integrability condition then
allows some $p>2$. A final disk Dirichlet estimate proves the assertion.
All exponent choices are strict, so no endpoint embedding is used.
\end{proof}

\begin{corollary}[boundary control of the Green function]\label{cor:greenC1}
For fixed $p\in\D$, the function
$U_p(z)=\cG(\Phi(p),\Phi(z))$ has bounded first derivatives in a neighborhood
of $\Sone$ and extends there with class $C^{1,\gamma}$ for some $\gamma>0$.
\end{corollary}

\begin{proof}
Choose a smooth cutoff $\chi$ equal to zero near $p$ and equal to one near
$\Sone$. The decomposition of $\cG$ in Lemma~\ref{lem:std} and conformal
invariance show that $v=\chi U_p\in H^1_0(\D)\cap L^\infty(\D)$:
the Laplace Green singularity has been removed, and the resolvent correction
is a bounded physical $H^1_0$ function. Away from $p$, $-\Delta U_p=aU_p$,
so $-\Delta v=av+f$ with
$f=-2\nabla\chi\cdot\nabla U_p-(\Delta\chi)U_p$ smooth and compactly
supported in $\D\setminus\{p\}$. Apply Lemma~\ref{lem:bootstrap}.
\end{proof}

The two preceding estimates do different jobs: Lemma~\ref{lem:gradient}
controls Neumann eigenfunctions in physical coordinates, whereas
Corollary~\ref{cor:greenC1} controls a Dirichlet Green function in disk
coordinates. Confusing these two boundary problems would not justify either use.

\begin{proposition}[after \cite{Nad88} and {\cite[Proposition~2.5]{OAI}}]\label{prop:dim}
$\dim V\le 2$.
\end{proposition}

\begin{proof}
We follow the boundary-sign proof of \cite[Lemmas~2.3, 2.4 and Proposition~2.5]{OAI} and indicate
the changes. Let $L\in H^1_0(\Omega)$ solve $(-\Delta-\mu)L=1$. Then $L\in H^2(\Omega)$ by
Lemma~\ref{lem:std}\ref{it:H2}, $L\ge0$ (test with the negative part and use $\mu<\lambda_D$),
$-\Delta L=1+\mu L>0$, so $L>0$ in $\Omega$, and $L$ is smooth up to the open smooth arcs. Hopf's lemma
gives $-\partial_\nu L>0$ on every open smooth arc. For $v\in V$, Green's identity for $H^2$ functions on
a Lipschitz domain gives
\[
0=\int_\Omega v=\int_\Omega v\,(-\Delta-\mu)L=-\int_{\partial\Omega}v\,\partial_\nu L\,d\sigma
=\int_{\Sone}\varpi\,(v\circ\Phi)\,ds,\qquad \varpi=|\Phi'|\,(-\partial_\nu L)\circ\Phi .
\]
Here $\varpi\in L^1(\Sone)$ and $\varpi>0$ on $A$. The trace $h_v=v\circ\Phi|_{\Sone}$ is continuous.
If $h_v=0$ then $v\in H^1_0(\Omega)$ and $\mu<\lambda_D$ forces $v=0$; otherwise $h_v$ takes both
signs. The sets $\{v>0\}$ and $\{v<0\}$ are connected. Indeed, if there were at
least three nodal domains, restrict $v$ to each of two of them and extend
these restrictions by zero to the rest of $\Omega$. These functions belong
to $H^1(\Omega)$ and retain Rayleigh quotient $\mu$, by truncation at the zero
level and the weak eigenfunction identity. They can be combined with nonzero
coefficients to have mean zero. Their Rayleigh quotient would still be $\mu$,
so the combination would be a $\mu$-eigenfunction vanishing on the third nodal
domain. Interior unique continuation rules this out. This is the usual
variational proof of the first-positive-eigenvalue case of Courant's argument. A boundary point where $v>0$ has a neighborhood in $\overline\Omega$ on which
$v>0$, and the local Lipschitz boundary chart joins it to nearby interior points
through the same sign region;
hence four cyclically ordered boundary points with strictly alternating signs would give two
disjoint crosscuts of $\overline\D$ (after applying the homeomorphism $\Phi^{-1}$) with alternating
endpoints, which is impossible. For completeness, we spell out the last step of the boundary-sign proof in
\cite[Proposition~2.5]{OAI}. Parametrize the circle by $s=e^{i\theta}$. If $\dim V\ge3$, choose $0\ne v\in V$ with
$\int\varpi h_v\cos\theta=\int\varpi h_v\sin\theta=0$ and write $h=h_v$.
Cut the circle at a strictly negative point of angle $\theta_0$. Let $\theta_-,\theta_+$ be the
infimum and supremum of the angles where $h>0$ in $(\theta_0,\theta_0+2\pi)$.
Continuity and the signs give $\theta_0<\theta_-<\theta_+<\theta_0+2\pi$. There is no positive
value outside $[\theta_-,\theta_+]$, and no negative value inside $(\theta_-,\theta_+)$: otherwise two positive
points on either side, together with the negative cut point, give four alternating signs.
Set
\[
F(\theta)=\cos\!\left(\theta-\frac{\theta_-+\theta_+}{2}\right)
-\cos\!\left(\frac{\theta_+-\theta_-}{2}\right).
\]
This function is positive on $(\theta_-,\theta_+)$ and negative outside $[\theta_-,\theta_+]$ in the cut
interval. Thus $Fh\ge0$ everywhere and $Fh>0$ on a nonempty open set. Since
$\varpi>0$ almost everywhere, $\int\varpi Fh>0$, contradicting the orthogonality
of $\varpi h$ to $1$, $\cos\theta$ and $\sin\theta$.
\end{proof}

\section{The vector principle and the curvature identity}\label{sec:rohleder}

A field $X\in H^1(\Omega;\R^2)$ is \emph{tangent} if the trace of $X\cdot\nu$ vanishes on every
open smooth arc. Write $\dv X=\partial_1X_1+\partial_2X_2$ and $\curl X=\partial_1X_2-\partial_2X_1$, and put $Q(X)=\int_\Omega((\dv X)^2+(\curl X)^2)$.

\begin{proposition}[Rohleder's principle {\cite{Roh21}}]\label{prop:rohleder}
For every tangent $X\in H^1(\Omega;\R^2)$,
\[
Q(X)\ \ge\ \mu\int_\Omega|X|^2,
\]
with equality if and only if $X=\nabla v$ for some $v\in V$.
\end{proposition}

\begin{proof}
Solve $\Delta\varphi_N=\dv X$ with $\partial_\nu\varphi_N=0$ and
$\int_\Omega\varphi_N=0$, and solve $\Delta\varphi_D=\curl X$ with
$\varphi_D\in H^1_0(\Omega)$. Tangency gives the Neumann compatibility
condition $\int_\Omega\dv X=0$. Both solutions are in $H^2(\Omega)$ by
Lemma~\ref{lem:std}\ref{it:H2}; for the Neumann solution apply that statement
to $(-\Delta+1)\varphi_N=-\dv X+\varphi_N\in L^2$.
For $\nabla^\perp\varphi_D=(-\partial_2\varphi_D,\partial_1\varphi_D)$,
the orientation convention gives
\[
(\nabla^\perp\varphi_D)\cdot\nu
=-\partial_\tau\varphi_D=0
\quad\text{on }\Gamma.
\]
Thus $Z=X-\nabla\varphi_N-\nabla^\perp\varphi_D$ is in $H^1$, tangent,
divergence free and curl free. On the simply connected $\Omega$ it has a
potential $Z=\nabla\varphi_H$, with $\varphi_H\in H^1(\Omega)$ after
adding a constant. The weak integration-by-parts identity gives
$\int_\Omega|\nabla\varphi_H|^2=0$. Hence
$X=\nabla\varphi_N+\nabla^\perp\varphi_D$. The two terms are
$L^2$-orthogonal, because $\dv(\nabla^\perp\varphi_D)=0$ and its normal
trace vanishes. Consequently
$Q(X)=\|\Delta\varphi_N\|_2^2+\|\Delta\varphi_D\|_2^2$.
Integration by parts, Cauchy--Schwarz and the Neumann Poincar\'e inequality give
\[
\|\nabla\varphi_N\|_2^2
=-\int_\Omega\varphi_N\Delta\varphi_N
\le\mu^{-1/2}\|\nabla\varphi_N\|_2\|\Delta\varphi_N\|_2.
\]
Together with the analogous Dirichlet inequality this yields
\[
Q(X)\ge\mu\|\nabla\varphi_N\|_2^2+
\lambda_D\|\nabla\varphi_D\|_2^2
=\mu\|X\|_2^2+(\lambda_D-\mu)\|\nabla\varphi_D\|_2^2.
\]
Equality forces $\varphi_D=0$ and equality in the Neumann Rayleigh
quotient, so $\varphi_N\in V$ (including $\varphi_N=0$).
Conversely, for $v\in V$, $\nabla v$ is an $H^1$ tangent field and
$Q(\nabla v)=\mu^2\|v\|_2^2=\mu\|\nabla v\|_2^2$.
\end{proof}

The spectral principle is expressed using divergence and curl. To transfer it
conformally we need its gradient form, including the signed curvature integral.
This is the curvature form underlying \cite[Section~3]{Roh21}; the curvature
convention there has the opposite sign to ours. We give the density argument
for all tangent $H^1$ fields explicitly.

\begin{lemma}[curvature identity without a vertex term]\label{lem:det}
For every tangent field $X\in H^1(\Omega;\R^2)$,
\begin{equation}\label{eq:Q}
Q(X)=\int_\Omega|\nabla X|^2+
\int_\Gamma\kappa|\operatorname{Tr}X|^2\,d\sigma.
\end{equation}
The trace integral is finite and there is no additional term supported at the vertices.
\end{lemma}

\begin{proof}
Pointwise, $(\dv X)^2+(\curl X)^2-|\nabla X|^2=2\det\nabla X$.
For a smooth tangent field that vanishes near the corners, the divergence
theorem and integration by parts along the boundary give
\[
2\int_\Omega\det\nabla X
=\int_\Gamma(X_1\partial_\tau X_2-X_2\partial_\tau X_1)\,d\sigma.
\]
On an arc write $X=f_\tau\tau$. The integrand is
$f_\tau^2\det(\tau,\partial_\tau\tau)=\kappa f_\tau^2$, with the sign convention fixed
above. This proves \eqref{eq:Q} for these fields, without endpoint terms.

We justify density in the tangent $H^1$ space. First truncate radially:
$X^\Lambda=X\min(1,\Lambda/|X|)$ for $X\ne0$, and $X^\Lambda=0$ where $X=0$. The radial retraction is Lipschitz,
$|\nabla X^\Lambda|\le|\nabla X|$, and $X^\Lambda\to X$ in $H^1$ as $\Lambda\to\infty$.
Its trace is still tangent. For fixed $\Lambda$, use smooth logarithmic cutoffs
$0\le\eta_\varepsilon\le1$ that vanish near the finite corner set, tend to
one almost everywhere, and satisfy $\|\nabla\eta_\varepsilon\|_2\to0$.
Such cutoffs follow explicitly from the planar logarithmic radial profile:
its squared Dirichlet norm between radii $\varepsilon^2$ and $\varepsilon$
is $2\pi/|\log\varepsilon|$, and smoothing preserves the vanishing limit.
Then
\[
\|\nabla(\eta_\varepsilon X^\Lambda-X^\Lambda)\|_2
\le\|(1-\eta_\varepsilon)\nabla X^\Lambda\|_2
+\Lambda\|\nabla\eta_\varepsilon\|_2\longrightarrow0,
\]
and convergence in $L^2$ is immediate.

A field vanishing near the corners can now be smoothed while preserving
tangency. Use a partition of unity by interior patches and tubular patches
of open arcs. In a tubular patch with coordinates $(s,t)$, express the field
as $f_\tau(s,t)\tau(s)+f_\nu(s,t)\nu(s)$. The normal coefficient has trace $f_\nu(s,0)=0$.
Extend $f_\tau$ evenly and $f_\nu$ oddly in $t$, mollify with an even kernel, and then
return to the smooth moving frame. The resulting field is smooth and tangent
on $t=0$, and converges in $H^1$ under the smooth coordinate change. Supports
can be kept inside the patches and away from the corners. Interior patches
use ordinary mollification. Their sum proves the required density.

Both sides of \eqref{eq:Q} are continuous under $H^1$ convergence: for the
boundary term use the boundedness of $\kappa$ and the continuous trace map
$H^1(\Omega)\to L^2(\partial\Omega)$. Thus the identity passes to every
tangent $H^1$ field.
\end{proof}

\begin{remark}
The turning-angle atoms of the geometric curvature measure are not inserted
in \eqref{eq:Q}. An arbitrary $H^1$ field has no canonical point values at
vertices, and the density proof determines the closed energy identity directly.
For eigenfunction gradients the continuous corner values are in fact zero,
but that stronger fact is not needed to prove this identity for all tangent fields.
\end{remark}

\section{Kernels for the singular conformal potential}\label{sec:kernels}

Keep the previously fixed $\Phi$ and $a$, and introduce the measure
\[
dq=a\,dx,\qquad a=\mu|\Phi'|^2 .
\]
By \eqref{eq:SC}, $a$ is locally bounded in $\D$ and near every point of $A$,
and has the integrable corner singularities in \eqref{eq:apotential}; $q(\D)=\mu|\Omega|<\infty$. Conformal invariance of the Dirichlet integral and the
Dirichlet inequality on $\Omega$ give
\begin{equation}\label{eq:sub}
\int_\D\psi^2\,dq\le d\int_\D|\nabla\psi|^2\qquad(\psi\in H^1_0(\D)),\qquad d=\mu/\lambda_D<1 .
\end{equation}
We write $T=T_q$, $K=K^q$, $R=R^q$, $N=N^q$ for the objects in \eqref{eq:kernels}, and
$\langle f,h\rangle_q=\int_\D fh\,dq$ and
$\|f\|_q=\langle f,f\rangle_q^{1/2}$. For any nonnegative measure $\rho$ used below, the analogous inner product
and norm on $L^2(\rho)$ use $d\rho$ instead of $dq$; for vector functions the integrand is
the Euclidean scalar product.

\begin{lemma}\label{lem:T}
For $f\in L^2(q)$ there is a unique $Tf\in H^1_0(\D)$ with
$\int_\D\nabla Tf\cdot\nabla\psi=\ip\psi f_q$ for all $\psi\in H^1_0(\D)$. The operator $T$ is
self-adjoint on $L^2(q)$ with $\|T\|\le d$, and $Tf(x)=\int_\D G(x,z)f(z)\,dq(z)$ for almost every
$x$; for $f\ge0$ this is the kernel integral of \eqref{eq:kernels}. The same holds for every
measure $0\le\rho\le q$, with the same $d$.
\end{lemma}

\begin{proof}
By \eqref{eq:sub}, $|\ip\psi f_q|\le\sqrt d\,\|\nabla\psi\|_2\|f\|_q$, so the Riesz theorem gives
$Tf$. Testing with $\psi=Tf$ and using \eqref{eq:sub} again,
$\|Tf\|_q^2\le d\|\nabla Tf\|_2^2=d\,\ip{Tf}f_q\le d\|Tf\|_q\|f\|_q$. Symmetry is clear from the weak
formulation. For $f\ge0$ bounded with compact support in $\D$ the source $af$ is bounded and the
Green representation is classical; a general $f\ge0$ is the increasing limit of such functions,
the potentials increase to the kernel integral, and they converge to $Tf$ in $H^1_0$.
\end{proof}

In particular, for each $0\le\rho\le q$ the bounded inverse on $L^2(\rho)$ is
\[
(I-T_\rho)^{-1}=\sum_{n=0}^{\infty}T_\rho^n,
\qquad \|(I-T_\rho)^{-1}\|\le(1-d)^{-1}.
\]
This follows from the preceding norm bound, before any solution formula uses
the inverse. The integral definition of $T_\rho$ initially applies to
nonnegative functions; the lemma identifies its bounded linear extension
to signed $L^2(\rho)$ data.

\subsection{The interior kernel is the Green function of \texorpdfstring{$\Omega$}{Omega}}
For $x\ne y$ in $\D$ let $M(x,y)=\sum_{n\ge0}(T^nG(\cdot,y))(x)\in(0,\infty]$, the sum over all
chains $G(x,z_1)G(z_1,z_2)\cdots G(z_n,y)$ integrated against $q^{\otimes n}$, and
\[
\Gq(x,y)=\cG(\Phi(x),\Phi(y)) .
\]

\begin{lemma}\label{lem:M}
$M=\Gq$ off the diagonal. In particular $M^\rho\le\Gq$ for every measure $0\le\rho\le q$, where
$M^\rho$ is the chain sum with $\rho$ in place of $q$.
\end{lemma}

\begin{proof}
The Green function of the Laplacian is conformally invariant,
$G^\Omega_0(\Phi(x),\Phi(z))=G(x,z)$, and $d\zeta=|\Phi'(z)|^2dz$. So \eqref{eq:resOmega} becomes
\begin{equation}\label{eq:res}
\Gq(x,y)=G(x,y)+\int_\D G(x,z)\,\Gq(z,y)\,dq(z)\qquad(x\ne y).
\end{equation}
For each integer $\ell\ge0$, $\ell+1$ substitutions of \eqref{eq:res} give
\[
\Gq(x,y)=\sum_{n=0}^{\ell}(T^nG(\cdot,y))(x)
+\ip{G(x,\cdot)}{T^{\ell}\Gq(\cdot,y)}_q.
\]
The remainder is nonnegative, so $M\le\Gq$. Both $G(x,\cdot)$ and
$\Gq(\cdot,y)$ belong to $L^2(q)$: each
has a logarithmic singularity at an interior point, where $a$ is bounded, the first is bounded
elsewhere, and $\|\Gq(\cdot,y)\|_q^2=\mu\|\cG(\cdot,\Phi(y))\|_{L^2(\Omega)}^2$. By
Lemma~\ref{lem:T} the remainder is at most $d^\ell\|G(x,\cdot)\|_q\|\Gq(\cdot,y)\|_q\to0$ as $\ell\to\infty$.
\end{proof}

\subsection{The boundary kernels}
Let $\cP(\xi,\eta)=-\partial_{\nu(\eta)}\cG(\xi,\eta)$ for
$\xi\in\Omega$, $\eta\in\Gamma$ be the Poisson kernel of $-\Delta-\mu$ on
$\Omega$. For $p\in\D$ and $s\ne t$ in $A$, put $\xi_s=\Phi(s)$,
$\xi_t=\Phi(t)$ and define
\[
\Kq(p,s)=|\Phi'(s)|\,\cP(\Phi(p),\xi_s),\qquad
\Nq(s,t)=|\Phi'(s)||\Phi'(t)|\,
\partial_{\nu(\xi_s)}\partial_{\nu(\xi_t)}\cG(\xi_s,\xi_t).
\]

\begin{lemma}\label{lem:bdry}
Let $p\in\D$ and $s\ne t$ in $A$.
\begin{enumerate}[label=\textup{(\alph*)},leftmargin=2em]
\item\label{it:K} $K(p,s)=\Kq(p,s)$; in particular $0<K(p,s)<\infty$.
\item\label{it:Kbdd} $\sup_{s\in A}K(p,s)<\infty$.
\item\label{it:N} $R(s,t)=\int_\D P_s(z)K(z,t)\,dq(z)$ and $N_0(s,t)\le N(s,t)\le\Nq(s,t)<\infty$.
\item\label{it:R} $R$ is symmetric and $\iint_{\Sone\times\Sone}R\,ds\,dt\le q(\D)/(1-d)$.
\end{enumerate}
\end{lemma}

\begin{proof}
\ref{it:K} By Tonelli and Lemma~\ref{lem:M},
\begin{equation}\label{eq:Ksum}
K(p,s)=P_s(p)+\int_\D M(p,z)P_s(z)\,dq(z)=P_s(p)+\int_\D\Gq(p,z)P_s(z)\,dq(z).
\end{equation}
Apply \eqref{eq:res}, with the roles of the variables exchanged by symmetry, at $y=rs$, $r\uparrow1$:
\[
\frac{\Gq(p,rs)}{1-r}=\frac{G(p,rs)}{1-r}+\int_\D\Gq(p,z)\,\frac{G(z,rs)}{1-r}\,dq(z).
\]
The explicit formula gives $G(z,rs)/(1-r)\to P_s(z)$ for every $z\in\D$. On the left, $\Phi$ is
smooth and conformal up to $s$, and $s\Phi'(s)=|\Phi'(s)|\nu(\xi_s)$.
Thus $\Phi(rs)=\xi_s-(1-r)|\Phi'(s)|\nu(\xi_s)+O((1-r)^2)$;
since $\cG(\Phi(p),\cdot)$ is smooth up to $\Gamma$ and vanishes there, the left side tends to
$\Kq(p,s)$. Fatou's lemma and \eqref{eq:Ksum} give $K(p,s)\le\Kq(p,s)<\infty$. For the equality
let $B_r$ be the disk with center $rs$ and radius $(1-r)/2$. If $z\in\D\setminus B_r$ then
$|z-s|\le3|z-rs|$ and
\[
\frac{G(z,rs)}{1-r}\le\frac{(1-|z|^2)(1+r)}{4\pi|z-rs|^2}\le 9P_s(z),
\]
and $\Gq(p,\cdot)P_s$ is $q$-integrable by what was just proved. If $z\in B_r$ then
$G(z,rs)\le\frac1{4\pi}\log(1+6(1-r)^2/|z-rs|^2)$, whose integral over $B_r$ is $O((1-r)^2)$, while
$a$ and $\Gq(p,\cdot)$ are bounded near $s$; this part contributes $O(1-r)$. Dominated convergence
gives $\Kq(p,s)=P_s(p)+\int\Gq(p,z)P_s(z)\,dq(z)=K(p,s)$.

\ref{it:Kbdd} By \ref{it:K}, on $A$ the kernel is the inward disk normal
derivative of $U_p(z)=\cG(\Phi(p),\Phi(z))$:
$K(p,s)=-\partial_{\nu_\D}U_p(s)$. Corollary~\ref{cor:greenC1} bounds these
derivatives uniformly, including as $s$ approaches a prevertex. This also
provides a continuous extension of $K(p,\cdot)|_A$ to the circle, although
only its boundedness on $A$ is used. No continuity or finiteness of $N$ at
the prevertices is asserted.

\ref{it:N} The first formula is Tonelli's theorem applied to \eqref{eq:kernels}. Apply
\eqref{eq:res} at $(rs,rt)$:
\[
\frac{\Gq(rs,rt)-G(rs,rt)}{(1-r)^2}=\int_\D\frac{G(rs,z)}{1-r}\cdot\frac{\Gq(z,rt)}{1-r}\,dq(z).
\]
By the computation in \ref{it:K}, the integrand tends to $P_s(z)\Kq(z,t)=P_s(z)K(z,t)$ for every
$z$. On the left, $G(rs,rt)/(1-r)^2\to N_0(s,t)$, and by the joint smoothness in
Lemma~\ref{lem:std}\ref{it:green} and the vanishing of $\cG$ on $\Gamma$ in each variable,
$\Gq(rs,rt)/(1-r)^2\to\Nq(s,t)$. Fatou's lemma gives $R\le\Nq-N_0$.

\ref{it:R} Symmetry follows by reversing the chains. Since $\int_{\Sone}P_s(z)\,ds=1$, Tonelli gives
$\iint R=\sum_n\ip1{T^n1}_q\le\sum_nd^n\|1\|_q^2$.
\end{proof}

\subsection{A global reciprocal-matrix flow}

We now prove the compact-support kernel assertion needed in the singular
limit. The initial regularization and the reciprocal-matrix strategy come
from \cite[Section~5]{OAI}; their proofs are included. The preservation step
is proved through a differential equation and invariance of the negative
semidefinite cone. This is an alternative proof of the finite resolvent assertion obtained by iterating the column-update
lemma \cite[Lemma~5.4]{OAI}. No assertion from
\cite[Theorem~5.1]{OAI} is assumed.

For real symmetric matrices, $J\succeq0$ and $J\preceq0$ mean positive and
negative semidefiniteness. Write $\mathsf P(J)$ when $J$ has at most one
positive eigenvalue, with multiplicity. This property is preserved by
congruence, principal submatrices, subtraction of a positive semidefinite
matrix, and limits in a fixed finite dimension. The notation
$\mathsf M^{\circ-1}$ means the entrywise reciprocal, not the matrix inverse.
Matrix norms without a subscript are Euclidean operator norms.
The entrywise product is denoted by $\odot$. If $U\succeq0$, then
$U\odot U\succeq0$: represent $U_{ij}$ by scalar products and take tensor
products of the representing vectors.

\begin{lemma}[global preservation]\label{lem:flow}
Let $\mathsf W=\operatorname{diag}(w_1,\ldots,w_n)$ with $w_i\ge0$.
Suppose $\mathsf{M}(\theta)$ is a real symmetric $C^1$ matrix on $[0,1]$, all its
entries are strictly positive, and
\[
\frac{d\mathsf{M}}{d\theta}=\mathsf{M}\mathsf W \mathsf{M}.
\]
If $\mathsf P(\mathsf{M}(0)^{\circ-1})$ holds, then
$\mathsf P(\mathsf{M}(\theta)^{\circ-1})$ holds for every $\theta\in[0,1]$.
\end{lemma}

\begin{proof}
Put $J=\mathsf{M}^{\circ-1}$; a dot denotes differentiation in $\theta$ in this proof.
Direct differentiation gives the locally Lipschitz vector field
\begin{equation}\label{eq:reciprocalflow}
\dot J_{ij}=-\sum_{k=1}^n w_k
 \frac{J_{ij}^2}{J_{ik}J_{jk}}
\end{equation}
on the open set of symmetric matrices with strictly positive entries.

First suppose $\mathsf P(J)$ holds, and take $v\in\ker J$.
For each $k$ let
\[
\mathsf D_k=\operatorname{diag}(J_{1k},\ldots,J_{nk}),\qquad
A_k=\mathsf D_k^{-1}J\mathsf D_k^{-1},\qquad
\sigma_k=J_{kk}^{-1}.
\]
The $k$-th row and column of $A_k$ are constant, equal to $\sigma_k>0$.
The Schur complement at that entry, together with $\mathsf P(A_k)$, gives
\[
U_k=\sigma_k\mathbf1\mathbf1^T-A_k\succeq0,
\]
where $\mathbf1$ is the vector of ones. More explicitly, the $k$-th row
of $U_k$ is zero and its remaining block is the negative of that Schur
complement. Set $y_k=\mathsf D_kv$. Then
\[
\mathbf1^Ty_k=(Jv)_k=0,\qquad A_ky_k=0,\qquad U_ky_k=0.
\]
Expanding $A_k\odot A_k$ and using these three identities yields
\begin{equation}\label{eq:nullflow}
v^T\dot Jv
=-\sum_k w_k\,y_k^T(U_k\odot U_k)y_k\le0.
\end{equation}
This identity applies at every matrix in the specified open set having
$\mathsf P$, not only at the initial matrix.

The passage from this sign condition at boundary points of the cone to
forward invariance is the tangency criterion for flow-invariant closed sets of
Nagumo \cite{Nag42} and Brezis \cite{Bre70}. We include the short argument for
the cone at hand. The case $n=1$ is immediate. For $n\ge2$,
fix the last index as a pivot and write
\[
J=\begin{pmatrix}J_{00}&b\\b^T&h\end{pmatrix},\qquad
h>0,\qquad S=J_{00}-bb^T/h.
\]
Here the index $00$ denotes the block of the first $n-1$ indices.
Congruence gives $J\sim\operatorname{diag}(S,h)$, so $\mathsf P(J)$ is
equivalent to $S\preceq0$. The variables $(S,b,h)$ are smooth coordinates
on the open set under consideration. Let $\mathcal F(S,b,h)$ be the
$S$ component of the vector field \eqref{eq:reciprocalflow}.
If $S\preceq0$ and $x\in\ker S$, the vector
$v=(x,-b^Tx/h)$ belongs to $\ker J$, and differentiating the Schur
complement gives
\[
x^T\mathcal F(S,b,h)x=v^T\dot Jv\le0.
\]

For a symmetric matrix $S$, let $S_+$ retain its positive eigenvalues
and replace the others by zero, and put $S_-=S-S_+$. Thus $S_-\preceq0$,
$S_+\succeq0$ and $S_-S_+=0$. Near any point with $S\preceq0$ in our
open coordinate domain, the projected point $(S_-,b,h)$ remains in that
domain. The vector field is Lipschitz on a sufficiently small compact
neighborhood containing these points. The squared distance to the negative
semidefinite cone is $\|S_+\|_F^2$, where
$\|B\|_F^2=\operatorname{tr}(B^TB)$ is the Frobenius norm. Its derivative is
\[
\frac12\frac{d}{d\theta}\|S_+\|_F^2
=\operatorname{tr}\bigl(S_+\mathcal F(S,b,h)\bigr).
\]
This derivative formula also holds at zero eigenvalues: the scalar
function $r\mapsto(\max(r,0))^2$ is $C^1$, or equivalently one may use
the differentiability of squared distance to a closed convex set.
The range of $S_+$ is contained in $\ker S_-$, so the null-vector
inequality at the projected point gives
$\operatorname{tr}(S_+\mathcal F(S_-,b,h))\le0$. Consequently, for a
local Lipschitz constant $L_*$,
\[
\frac12\frac{d}{d\theta}\|S_+\|_F^2
\le L_*\|S_+\|_F^2.
\]
Gronwall's inequality shows that this distance remains zero if it starts
at zero. This proves local forward invariance at every cone point in the
open domain. A first-exit argument extends it along the whole trajectory:
at any alleged first exit the trajectory is still in the open domain and
the same local argument applies. Positivity and continuity of the entries
on the given compact time interval ensure that this domain is never left.
Thus $S\preceq0$ throughout, proving the lemma.
\end{proof}

\begin{lemma}[a finite resolvent]\label{lem:finiteres}
Let $\mathsf C$ be real symmetric with strictly positive entries,
$\mathsf P(\mathsf C^{\circ-1})$, and let $\mathsf W$ be diagonal with
nonnegative entries. If
$\|\mathsf W^{1/2}\mathsf C\mathsf W^{1/2}\|<1$, then the matrix
\[
\mathsf C+\mathsf C\mathsf W^{1/2}
(I-\mathsf W^{1/2}\mathsf C\mathsf W^{1/2})^{-1}
\mathsf W^{1/2}\mathsf C
\]
has an entrywise reciprocal with property $\mathsf P$.
\end{lemma}
\begin{proof}
Put $\mathsf B=\mathsf W^{1/2}\mathsf C\mathsf W^{1/2}$ and
\[
\mathsf{M}(\theta)=\mathsf C+\theta\mathsf C\mathsf W^{1/2}
(I-\theta\mathsf B)^{-1}\mathsf W^{1/2}\mathsf C,\qquad0\le\theta\le1.
\]
The inverse exists on this whole interval. Its geometric series shows
that $\mathsf{M}$ has strictly positive entries, since $\mathsf B$ has nonnegative
entries. Differentiation, or the identity
$\mathsf{M}(\theta)\mathsf W^{1/2}
=\mathsf C\mathsf W^{1/2}(I-\theta\mathsf B)^{-1}$, gives
$\dot{\mathsf M}=\mathsf{M}\mathsf W \mathsf{M}$. Apply Lemma~\ref{lem:flow} at $\theta=1$.
Neither invertibility nor positive definiteness of $\mathsf C$ is needed.
\end{proof}

\subsection{An initial regularization, with its proof}

The following construction is adapted from \cite[Section~5]{OAI}. It is
reproved to avoid taking its reciprocal-kernel assertion as an input.

\begin{lemma}[regularization from \cite{OAI}]\label{lem:initreg}
There are bounded continuous strictly positive kernels $C_\eta$ on
$\D\times\D$, $0<\eta\le1$, such that $C_\eta\uparrow G$ as $\eta\downarrow0$,
including the infinite diagonal, and every finite entrywise reciprocal
matrix of $C_\eta$ has property $\mathsf P$.
\end{lemma}
\begin{proof}
First define
\[
\mathcal L(r)=\frac1{\log(1+1/r)}\quad(r>0),\qquad\mathcal L(0)=0.
\]
For $x,y\in\D$ and $r\ge0$ put
\[
\mathfrak h(x,y)=\frac{|x-y|^2}{(1-|x|^2)(1-|y|^2)},\qquad
f_*(r)=r+\tfrac12-\mathcal L(r).
\]
Positive diagonal congruence with entries $1-|x|^2$ turns the kernel
$\mathfrak h+\tfrac12$ into
$\tfrac12(1+|x|^2)(1+|y|^2)-2x\cdot y$.
This is a positive rank-one kernel minus a positive semidefinite kernel
and hence has $\mathsf P$ on every finite set.

We next show that $f_*(\mathfrak h)$ is positive semidefinite. For $0<\alpha<1$,
normalize $v^{\alpha-1}(1+v)^{-1}\,dv$ to a probability measure on $(0,\infty)$
and denote expectation by $\mathbb E_\alpha$. Write $r_v=(1+v)^{-1}$.
Integration by parts gives $\mathbb E_\alpha(1-r_v)=\alpha$.
The elementary identities
\[
\mathcal L(r)=\int_0^1 r^{1-\alpha}(r+1)^\alpha\,d\alpha,\qquad
r^{1-\alpha}(r+1)^\alpha
=\mathbb E_\alpha\frac{r(r+1)}{r+r_v}
\]
hold for $r>0$. The second follows by scaling $v$ in the defining
integral: the normalization cancels
$\int_0^\infty v^{\alpha-1}/(1+v)\,dv$. Polynomial division then gives
\[
f_*(r)=\int_0^1\mathbb E_\alpha
\frac{r_v(1-r_v)}{r+r_v}\,d\alpha,\qquad
0\le f_*(r)\le f_*(0)=\tfrac12.
\]
The formula extends to $r=0$ by dominated convergence. Now
\[
k_*(x,y)=\frac1{1+\mathfrak h(x,y)}
=\frac{(1-|x|^2)(1-|y|^2)}{|1-x\overline y|^2}
\]
is positive semidefinite: the geometric power series makes
$(1-x\overline y)^{-1}$ a complex Gram kernel, and multiplication by its
complex conjugate and the diagonal factors preserves that property.
For every $0<r_v<1$,
\[
\frac1{\mathfrak h+r_v}
=\sum_{j=0}^\infty(1-r_v)^j k_*^{\,j+1}
\]
is a sum of positive semidefinite kernels, where powers are entrywise.
It converges also on the diagonal. The positive mixture above proves
positive semidefiniteness of $f_*(\mathfrak h)$.

Define
\[
\frac1{C_\eta(x,y)}
=4\pi\bigl[\mathfrak h+\tfrac12-(1-\eta)f_*(\mathfrak h)\bigr]
=4\pi\bigl[\mathcal L(\mathfrak h)+\eta f_*(\mathfrak h)\bigr].
\]
The first expression has $\mathsf P$. The bracket is at least $\eta/2$,
since it equals
$\eta(\mathfrak h+\tfrac12)+(1-\eta)\mathcal L(\mathfrak h)$.
Thus $C_\eta$ is bounded, continuous and strictly positive.
Since $f_*\ge0$, it increases to $G$ as $\eta\downarrow0$, and
$C_\eta(x,x)=1/(2\pi\eta)$.
\end{proof}

\subsection{From the flow to the compact-support boundary theorem}

\begin{proposition}[compact-support reciprocal theorem]\label{prop:compactkernel}
Let $d\rho=a_\rho\,dx$, where $a_\rho\ge0$ is bounded and supported in a
compact subset of $\D$. Suppose
\[
\int_\D\psi^2\,d\rho\le d_\rho\int_\D|\nabla\psi|^2
\quad(\psi\in H^1_0(\D)),\qquad 0<d_\rho<1.
\]
For $p\in\D$ define
\[
D_p^\rho(s,t)=\frac{K^\rho(p,s)K^\rho(p,t)}{N^\rho(s,t)}\quad(s\ne t),
\qquad D_p^\rho(s,s)=0.
\]
Then $D_p^\rho$ is conditionally negative semidefinite on $\Sone$.
\end{proposition}

\begin{proof}
The proof of Lemma~\ref{lem:T}, with this form inequality in place of
\eqref{eq:sub}, gives a self-adjoint $T_\rho$ on $L^2(\rho)$ with
$\|T_\rho\|\le d_\rho<1$. For distinct interior points define
\[
M^\rho(x,y)=G(x,y)+
\langle G(x,\cdot),(I-T_\rho)^{-1}G(y,\cdot)\rangle_\rho.
\]
Its correction is finite: logarithmic singularities belong to $L^2(\rho)$.
Set $M^\rho(x,x)=+\infty$.
We first prove $\mathsf P$ for its finite reciprocal matrices with zero
diagonal.

Fix $0<\eta\le1$ and replace $G$ by $C_\eta$ from Lemma~\ref{lem:initreg}.
The resulting integral operator $T_\eta$ satisfies
$|T_\eta f|\le T_\rho|f|$, and hence $\|T_\eta\|\le d_\rho$.
This estimate does not assume positive semidefiniteness of $T_\eta$.
Choose a compact integration set $Q_\rho\Subset\D$ containing the support
and partition it into finitely many measurable cells
whose maximal diameter tends to zero. Choose a representative $z_\alpha$ in every nonempty cell and let
$w_\alpha$ be its mass. Cells of zero $\rho$-mass are omitted only from
the finite matrices and normalized indicators below. Let $\Pi_{\mathrm{cell}}$ map each cell to its representative.
Uniform continuity on this compact set gives
\[
\|C_\eta(\Pi_{\mathrm{cell}}(\cdot),\Pi_{\mathrm{cell}}(\cdot))-C_\eta\|_\infty\longrightarrow0.
\]
The corresponding operators on $L^2(\rho)$ therefore converge to $T_\eta$
in norm, with error bounded by $\rho(\D)$ times that uniform error.
Their restrictions to the span of the normalized cell indicators are the symmetric
matrices
$\mathsf B_{\alpha\beta}
=\sqrt{w_\alpha}C_\eta(z_\alpha,z_\beta)\sqrt{w_\beta}$;
they vanish on the orthogonal complement. In particular $\|\mathsf B\|<1$
for all sufficiently fine partitions. More precisely, their operator norms
are eventually at most $(1+d_\rho)/2$, so their resolvent norms are
at most $2/(1-d_\rho)$.

Fix any finite set of evaluation points and join it to the integration
nodes, merging coincident points. Give evaluation-only points weight zero.
The matrix $\mathsf C$ of $C_\eta$ on this union has strictly positive
entries and $\mathsf P(\mathsf C^{\circ-1})$ by Lemma~\ref{lem:initreg}.
After removing the rows and columns belonging to evaluation-only points,
the remaining block of $\mathsf W^{1/2}\mathsf C\mathsf W^{1/2}$ is
$\mathsf B$; the removed rows and columns are zero. Lemma~\ref{lem:finiteres} therefore proves
$\mathsf P$ for the reciprocal of the discrete resolvent, on the whole
union and then on the evaluation points.

For fixed $x$, the source $C_\eta(x,\Pi_{\mathrm{cell}}(\cdot))$ converges uniformly on
the integration set, hence in $L^2(\rho)$, to $C_\eta(x,\cdot)$.
Norm convergence of the operators, their common strict norm bound for
fine partitions, and the identity for the difference of inverses show
that the discrete resolvent entries converge to
\[
M^\rho_\eta(x,y)=C_\eta(x,y)+
\langle C_\eta(x,\cdot),(I-T_\eta)^{-1}C_\eta(y,\cdot)\rangle_\rho.
\]
Thus its reciprocal matrices have $\mathsf P$. These partitions are an
analytic limiting argument; no finite computed partition is used as evidence.

Let $\eta\downarrow0$. The geometric resolvent series is the sum of
nonnegative chain integrals. Monotone convergence in each product integral
and then in the series gives $M^\rho_\eta(x,y)\uparrow M^\rho(x,y)$ for
$x\ne y$. The correction in the limit is bounded in absolute value by
\[
\frac{\|G(x,\cdot)\|_{L^2(\rho)}
       \|G(y,\cdot)\|_{L^2(\rho)}}{1-d_\rho}.
\]
On the diagonal the direct term tends to infinity. Taking reciprocal
entries and fixed-size limits proves $\mathsf P$ for $M^\rho$ with
reciprocal diagonal zero. If $\rho=0$, the same conclusion follows
directly from the regularization without partitions.

On the compact integration support the Poisson sections are bounded, so
Tonelli's theorem and the resolvent series give
\[
\begin{aligned}
K^\rho(p,s)&=P_s(p)+\langle G(p,\cdot),(I-T_\rho)^{-1}P_s\rangle_\rho,\\
N^\rho(s,t)&=N_0(s,t)+\langle P_s,(I-T_\rho)^{-1}P_t\rangle_\rho
\quad(s\ne t).
\end{aligned}
\]
Now fix an integer $n_b\ge1$, take distinct $t_1,\ldots,t_{n_b}\in\Sone$, keep $p\in\D$ fixed, and let
$r\uparrow1$ with $r>|p|$. The explicit Green formula gives
\[
\frac{G(rs,x)}{1-r}\longrightarrow P_s(x),\qquad
\frac{G(rs,rt)}{(1-r)^2}\longrightarrow N_0(s,t)\quad(s\ne t).
\]
The first convergence is uniform on the compact integration support.
Using the bounded resolvent and the $L^2(\rho)$ source at $p$ gives
\[
\frac{M^\rho(p,rt_i)}{1-r}\longrightarrow K^\rho(p,t_i),\qquad
\frac{M^\rho(rt_i,rt_j)}{(1-r)^2}\longrightarrow N^\rho(t_i,t_j)
\quad(i\ne j).
\]
All limiting entries are finite and positive: the resolvent expressions
use bounded Poisson sections on the support. Congruence of the reciprocal
matrix by diagonal factors $1-r$ at moving points and $1$ at $p$, followed
by the limit, gives $\mathsf P$ for
\[
\mathsf J=\begin{pmatrix}B&\ell\\\ell^T&0\end{pmatrix},\qquad
B_{ii}=0,\quad B_{ij}=\frac1{N^\rho(t_i,t_j)}\ (i\ne j),\quad
\ell_i=\frac1{K^\rho(p,t_i)}.
\]

For completeness, this implies $x^TBx\le0$ when $\ell^Tx=0$.
Otherwise choose $y$ with $\ell^Ty\ne0$ and subtract
$(x^TBy)/(x^TBx)$ times $x$ to arrange $x^TBy=0$.
Then $(x,0)$ and $(y,c)$ are orthogonal for the form $\mathsf J$.
Choosing $c$ so that $y^TBy+2c\ell^Ty>0$ produces a two-dimensional
positive subspace, contrary to $\mathsf P$.
Finally set $x_i=K^\rho(p,t_i)\omega_i$ for coefficients with
$\sum_i\omega_i=0$. The resulting inequality is exactly the required
conditional negativity. Repeated boundary points are handled by combining
their coefficients.
\end{proof}

The proof controls every fixed finite configuration at full subcritical
amplitude: the only strict bound used is the operator bound $d_\rho<1$, before
the flow and in the operator limits.

\subsection{Conditional negativity for the singular potential}

\begin{proposition}\label{prop:cnd}
For every $p\in\D$ the kernel
\[
D_p(s,t)=\frac{K(p,s)K(p,t)}{N(s,t)}\quad(s\ne t),\qquad D_p(s,s)=0,
\]
is conditionally negative semidefinite on $A$. Moreover, with $C_p=\sup_AK(p,\cdot)$,
\begin{equation}\label{eq:Dbound}
0<D_p(s,t)\le\pi C_p^2|s-t|^2\qquad(s\ne t\ \text{in}\ A).
\end{equation}
\end{proposition}

\begin{proof}
Let $r_n\uparrow1$ and $\rho_n=\mathbf 1_{\{|x|\le r_n\}}q$. Each $\rho_n$ has a bounded density and
satisfies \eqref{eq:sub} with the same $d$. The finite-valued assertions for these
truncations require no additional theorem from \cite{OAI}: on their compact
support the Poisson functions $P_s$ are uniformly bounded for $s\in\Sone$.
The formulas
\[
K^{\rho_n}(p,s)=P_s(p)+
\langle G(p,\cdot),(I-T_{\rho_n})^{-1}P_s\rangle_{\rho_n},
\qquad
R^{\rho_n}(s,t)=\langle P_s,(I-T_{\rho_n})^{-1}P_t\rangle_{\rho_n}
\]
are finite by Lemma~\ref{lem:T}; $G(p,\cdot)\in L^2(\rho_n)$ even if the
pole lies in the support. Thus Proposition~\ref{prop:compactkernel} applies and gives conditional negative semidefiniteness of $D^{\rho_n}_p$ on $\Sone$. Every term in \eqref{eq:kernels} is an integral of a nonnegative
function against a power of the measure, so by monotone convergence $K^{\rho_n}(p,s)\uparrow K(p,s)$
and $N^{\rho_n}(s,t)\uparrow N(s,t)$ for all $s\ne t$. On $A$ the limits are finite and positive by
Lemma~\ref{lem:bdry}, hence $D^{\rho_n}_p(t_i,t_j)\to D_p(t_i,t_j)$ for any finite family $(t_i)$ in $A$,
and the inequality passes to the limit. The bound \eqref{eq:Dbound} follows from $N\ge N_0$ and
Lemma~\ref{lem:bdry}\ref{it:Kbdd}.
\end{proof}

\subsection{Solution formulas}
For a bounded Borel function $h$ on $\Sone$ let $\mathsf{H}h(x)=\int_{\Sone}P_s(x)h(s)\,ds$ be its bounded
harmonic extension. Since $q$ is finite, $\mathsf{H}h\in L^2(q)$, and we define
\begin{equation}\label{eq:W}
W_h=\mathsf{H}h+T(I-T)^{-1}\mathsf{H}h .
\end{equation}

\begin{lemma}\label{lem:W}
Let $h,k$ be bounded Borel functions on $\Sone$.
\begin{enumerate}[label=\textup{(\alph*)},leftmargin=2em]
\item\label{it:Weq} $\zeta_h=W_h-\mathsf{H}h\in H^1_0(\D)$ and $\int_\D\nabla \zeta_h\cdot\nabla\psi=\ip\psi{W_h}_q$
for all $\psi\in H^1_0(\D)$.
\item\label{it:Wrep} $W_h(p)=\int_{\Sone}K(p,s)h(s)\,ds$ for every $p\in\D$ (continuous
representative).
\item\label{it:Rrep} $\ip{\mathsf{H}h}{(I-T)^{-1}\mathsf{H}k}_q=\iint_{\Sone\times\Sone}R(s,t)h(s)k(t)\,ds\,dt$.
\end{enumerate}
\end{lemma}

\begin{proof}
\ref{it:Weq} holds because $W_h=(I-T)^{-1}\mathsf{H}h$ in $L^2(q)$ and $\zeta_h=T(I-T)^{-1}\mathsf{H}h$. For \ref{it:Wrep}
and \ref{it:Rrep} expand $(I-T)^{-1}=\sum_nT^n$, which converges in operator norm, write
$T^n\mathsf{H}h=\int(T^nP_s)h(s)\,ds$ and exchange sums and integrals. For $|h|,|k|$ all terms are
nonnegative and the totals are bounded by $\|h\|_\infty\int K(p,s)\,ds$ and
$\|h\|_\infty\|k\|_\infty\iint R$, finite by Lemma~\ref{lem:bdry}; this justifies the exchanges for
signed data. Pointwise evaluation at $p$ uses $\zeta_h(p)=\ip{G(p,\cdot)}{(I-T)^{-1}\mathsf{H}h}_q$. Here
$p\mapsto G(p,\cdot)$ is continuous from $\D$ into $L^2(q)$ (the potential is bounded near $p$, and
locally the logarithmic singularity is continuous under translation in
$L^2$, and away from $p$ dominated convergence applies because the
sections are uniformly bounded while $q$ is finite), so this formula is the
continuous representative of $\zeta_h$ and \ref{it:Wrep} holds at every point, not only almost
everywhere.
\end{proof}

\section{Boundary energy and the multiplier identity}\label{sec:energy}

\begin{definition}
$\cA$ is the vector space of bounded Borel maps $h:\Sone\to\R^2$ that belong to
$H^{1/2}(\Sone;\R^2)$ and are \emph{tangent}: $h(s)\cdot\nu(\Phi(s))=0$ for almost every $s\in A$. For
$h\in\cA$ let $W_h$ be defined componentwise by \eqref{eq:W} and
\[
X_h=W_h\circ\Phi^{-1}.
\]
Only the almost-everywhere class of $h$ affects $W_h$, $X_h$, or the energy
below. All pointwise claims about eigenfunction traces will use their
continuous representatives. The next lemma establishes Sobolev regularity
before an energy is assigned to $X_h$.
\end{definition}

\begin{lemma}[well-defined physical extension]\label{lem:X}
For every $h\in\cA$, the field $X_h$ belongs to $H^1(\Omega;\R^2)$ and is
tangent, with trace $h\circ\Phi^{-1}$ on $\Gamma$.
\end{lemma}
\begin{proof}
Since $h\in H^{1/2}$, $\mathsf{H}h\in H^1(\D)$ with trace $h$, and $\zeta_h\in H^1_0(\D)$, so
$W_h\in H^1(\D)$. The Dirichlet integral is conformally invariant and
$\int_\Omega|X_h|^2=\mu^{-1}\|W_h\|_q^2<\infty$, so $X_h\in H^1(\Omega)$. Near a compact subset of
an open smooth arc $\Phi^{-1}$ is a diffeomorphism up to the boundary, so traces correspond: the trace of
$X_h$ on $\Gamma$ is $h\circ\Phi^{-1}$, which is tangent.

\end{proof}

\begin{definition}\label{def:energy}
For $h\in\cA$, define
\[
E(h)=Q(X_h)-\mu\|X_h\|_{L^2(\Omega)}^2,
\qquad E(h,k)=\tfrac12\big(E(h+k)-E(h)-E(k)\big).
\]
These quantities are finite by Lemma~\ref{lem:X}. The linearity of
$h\mapsto X_h$ makes $E(h,k)$ a symmetric bilinear form, with $E(h,h)=E(h)$.
\end{definition}

In \cite{OAI} the data are Lipschitz. The corner gradient here is continuous
and bounded, but Lipschitz regularity is not required. Its Sobolev trace is
$H^{1/2}$, which motivates this class. Define also the pulled-back curvature
\[
c(s)=\kappa(\Phi(s))|\Phi'(s)|\quad(s\in A),
\qquad \int_{\Sone}|c|\,ds=\int_\Gamma|\kappa|\,d\sigma<\infty.
\]
Values assigned at the finite exceptional set do not affect any integral.
Although $\kappa$ is bounded on the physical arcs, $c$ need not be bounded.

\begin{lemma}\label{lem:E}
Let $h,k\in\cA$.
\begin{enumerate}[label=\textup{(\alph*)},leftmargin=2em]
\item\label{it:pos} $E(h)\ge0$, and $E(h)=0$ if and only if $X_h=\nabla v$ for some $v\in V$.
\item\label{it:form} One has
\begin{equation}\label{eq:E}
\begin{aligned}
E(h,k)={}&\frac12\iint N_0(s,t)(h(s)-h(t))\cdot(k(s)-k(t))\,ds\,dt\\
&-\iint R(s,t)h(s)\cdot k(t)\,ds\,dt
+\int_{\Sone}c(s)h(s)\cdot k(s)\,ds .
\end{aligned}
\end{equation}
\end{enumerate}
\end{lemma}

\begin{proof}
\ref{it:pos} is Proposition~\ref{prop:rohleder}.

\ref{it:form} By \eqref{eq:Q} and conformal invariance,
$E(h,k)=\int_\D\nabla W_h:\nabla W_k-\ip{W_h}{W_k}_q
+\int_{\Sone}c\,h\cdot k$. The last term is absolutely integrable because
$h,k$ are bounded and $c\in L^1$. For the volume terms, by Lemma~\ref{lem:W}\ref{it:Weq} for $k$,
tested with $\zeta_h$, the function $W_h$ may be replaced by $\mathsf{H}h$; since $\mathsf{H}h$ is harmonic and
$\zeta_k\in H^1_0$, $\int\nabla \mathsf{H}h:\nabla W_k=\int\nabla \mathsf{H}h:\nabla \mathsf{H}k$. Thus
\[
E(h,k)=\int_\D\nabla \mathsf{H}h:\nabla \mathsf{H}k-\ip{\mathsf{H}h}{(I-T)^{-1}\mathsf{H}k}_q+\int_{\Sone}c\,h\cdot k .
\]
The first term is the Douglas integral, valid for $H^{1/2}$ data (both sides diagonalize in Fourier
modes), and the second is Lemma~\ref{lem:W}\ref{it:Rrep}. The first double
integral is absolutely convergent by Cauchy--Schwarz in the $H^{1/2}$
seminorm; the second is absolutely convergent because $R\in L^1$ and
$h,k$ are bounded. The last term is precisely the transformed curvature integral.
\end{proof}

Fix $0\ne u\in V$ and put $g=(\nabla u)\circ\Phi|_{\Sone}$.

\begin{lemma}\label{lem:g}
$g\in\cA$, $g$ is continuous, $g(s_j)=0$ for every prevertex, $W_g=(\nabla u)\circ\Phi$, and
\begin{equation}\label{eq:rep}
\int_{\Sone}K(p,s)\,g(s)\,ds=\nabla u(\Phi(p))\qquad(p\in\D),\qquad E(g)=0 .
\end{equation}
If $b:\Sone\to\R$ is Lipschitz, then $bg\in\cA$.
\end{lemma}

\begin{proof}
By Lemma~\ref{lem:gradient} each $w=\partial_iu$ is in $C(\overline\Omega)\cap H^1(\Omega)$,
vanishes at the vertices and solves $(\Delta+\mu)w=0$; $\nabla u$ is tangent on the open smooth arcs.
Then $w\circ\Phi\in C(\overline\D)\cap H^1(\D)$, so its restriction $g_i$ to $\Sone$ is in
$H^{1/2}$; hence $g\in\cA$. For $\psi\in C_c^\infty(\D)$, a change of variables gives
$\int_\D\nabla(w\circ\Phi)\cdot\nabla\psi=\ip\psi{w\circ\Phi}_q$, and this extends to
$\psi\in H^1_0(\D)$ by \eqref{eq:sub}. By Lemma~\ref{lem:W}\ref{it:Weq}, $W_{g_i}$ satisfies the
same identity and has the same trace. The difference $\delta\in H^1_0(\D)$ satisfies
$\int|\nabla\delta|^2=\|\delta\|_q^2\le d\int|\nabla\delta|^2$, so $\delta=0$. Both $W_g$ (in the
representative of Lemma~\ref{lem:W}\ref{it:Wrep}) and $(\nabla u)\circ\Phi$ are continuous in $\D$,
so they agree at every $p\in\D$; this is the first part of \eqref{eq:rep}. The second part follows
from $X_g=\nabla u$ and Proposition~\ref{prop:rohleder}. Finally, if $b$ has Lipschitz constant $\ell$, then
$|b(s)g(s)-b(t)g(t)|^2\le2\|b\|_\infty^2|g(s)-g(t)|^2+2\ell^2\|g\|_\infty^2|s-t|^2$, and
$N_0(s,t)|s-t|^2$ is bounded, so $bg\in H^{1/2}$; it is bounded and tangent.
\end{proof}

\begin{proposition}[multiplier identity]\label{prop:mult}
For every Lipschitz $b:\Sone\to\R$,
\begin{equation}\label{eq:mult}
E(bg)=\frac12\iint_{\Sone\times\Sone}N(s,t)\,(b(s)-b(t))^2\,g(s)\cdot g(t)\,ds\,dt ,
\end{equation}
and the integral converges absolutely.
\end{proposition}

\begin{proof}
$E$ is a nonnegative quadratic form on $\cA$ and $E(g)=0$, hence $E(g,k)=0$ for all $k\in\cA$; in
particular $E(g,b^2g)=0$. Subtract \eqref{eq:E} for the pair $(g,b^2g)$ from \eqref{eq:E} for
$(bg,bg)$. The curvature contribution cancels exactly:
\begin{equation}\label{eq:cancel}
\int_{\Sone}c(s)|b(s)g(s)|^2\,ds
-\int_{\Sone}c(s)g(s)\cdot b(s)^2g(s)\,ds=0.
\end{equation}
Both integrals are absolutely convergent, and this uses no sign condition on
$c$ or $\kappa$. Both $N_0$-integrals converge absolutely (the data are in $H^{1/2}$), and the difference
of the integrands, writing $b_s=b(s)$, $b_t=b(t)$, $g_s=g(s)$ and $g_t=g(t)$, is
\[
|b_sg_s-b_tg_t|^2-(g_s-g_t)\cdot(b_s^2g_s-b_t^2g_t)=(b_s-b_t)^2\,g_s\cdot g_t .
\]
For the second terms, the symmetry of $R$ gives
$-\iint R\,(b_sb_t-b_s^2)\,g_s\cdot g_t=\frac12\iint R\,(b_s-b_t)^2g_s\cdot g_t$. Absolute
convergence: $N_0(b_s-b_t)^2$ is bounded and $R$ is integrable.
\end{proof}

\section{Proof of Theorem A}\label{sec:proof}

\begin{lemma}\label{lem:embed}
Let $p\in\D$. There are a separable real Hilbert space $\cH$ and a Lipschitz map
$\overline{\mathcal B}:\Sone\to\cH$, injective on $A$, with $\|\overline{\mathcal B}(s)-\overline{\mathcal B}(t)\|^2=D_p(s,t)$ for $s,t\in A$.
\end{lemma}

\begin{proof}
We recall the Gram-kernel construction in Schoenberg's theorem \cite{Sch38}.
Fix $s_*\in A$ and set
\[
H_p(s,t)=\tfrac12\bigl(D_p(s,s_*)+D_p(t,s_*)-D_p(s,t)\bigr).
\]
For any finite points $t_i\in A$ and coefficients $c_i$, append the coefficient $-\sum_i c_i$
at $s_*$ in Proposition~\ref{prop:cnd}. The resulting inequality is
$\sum_{i,j}c_ic_jH_p(t_i,t_j)\ge0$. Thus $H_p$ is a positive semidefinite
Gram kernel. Quotient the real vector space spanned by symbols $e_s$
by the null space of the bilinear form $\langle e_s,e_t\rangle=H_p(s,t)$,
and complete it to obtain $\cH$. The images $\mathcal B(s)$ of $e_s$ satisfy
$\|\mathcal B(s)-\mathcal B(t)\|^2=D_p(s,t)$. By
\eqref{eq:Dbound} the resulting map $\mathcal B$ is injective and Lipschitz on $A$. Since $A$ is dense in $\Sone$, it extends
uniquely to a Lipschitz map $\overline{\mathcal B}$ on $\Sone$. The image
is compact and hence separable; its closed linear span is a separable Hilbert
space and contains the whole construction. Restrict $\cH$ to this span.
\end{proof}

The extension need not be injective at the prevertices. The next elementary lemma replaces the
statement ``an injective continuous image of a circle is not contained in a line''.

\begin{lemma}\label{lem:top}
No continuous $f:\Sone\to\R$ is injective on the complement of a finite set $F$.
\end{lemma}

\begin{proof}
Such an $f$ is not constant; let $f(s_-)=\min f<\max f=f(s_+)$. Each of the two open arcs of
$\Sone\setminus\{s_-,s_+\}$ is mapped onto a set containing $(\min f,\max f)$. Choose $t_0$ in this
interval outside the finite set $f(F)$. Then $t_0$ is attained in both arcs, at two distinct points
outside $F$.
\end{proof}

\begin{proof}[Proof of Theorem~\ref{thm:A}]
Let $0\ne u\in V$, $g$ as above, and suppose $\nabla u(\Phi(p))=0$ for some $p\in\D$. Take
$\overline{\mathcal B}$ from Lemma~\ref{lem:embed} and let $b_j$ be its coordinates in an orthonormal basis of
$\cH$; they are Lipschitz. By Proposition~\ref{prop:mult}, for partial sums the absolute integrand
is at most
\[
N(s,t)\,D_p(s,t)\,|g(s)||g(t)|=K(p,s)K(p,t)|g(s)||g(t)|\le C_p^2\|g\|_\infty^2
\]
for almost every $(s,t)$: the diagonal and pairs involving a prevertex
form an arclength-product null set. Thus by dominated convergence and \eqref{eq:rep}
\[
\sum_jE(b_jg)=\frac12\iint K(p,s)K(p,t)\,g(s)\cdot g(t)\,ds\,dt
=\frac12\Big|\int_{\Sone}K(p,s)g(s)\,ds\Big|^2=\frac12|\nabla u(\Phi(p))|^2=0 .
\]
Every term is nonnegative, so $E(b_jg)=0$ for all $j$, and by Lemma~\ref{lem:E} there are
$u_j\in V$ with $X_{b_jg}=\nabla u_j$. The traces agree almost everywhere on $\partial\Omega$ and
the boundary functions are continuous. On every compact subarc of $A$,
$\Phi$ is a smooth boundary diffeomorphism, so equality holds everywhere
there. Continuity extends it to the prevertices as well; both sides vanish
there by Lemma~\ref{lem:gradient}. Thus
\begin{equation}\label{eq:traces}
(\nabla u_j)\circ\Phi=b_j\,g\qquad\text{on }\Sone .
\end{equation}
Let $S_g=\{s\in\Sone:g(s)\ne0\}$. It is open, and contained in $A$. It is nonempty: otherwise each
$\partial_iu\in C(\overline\Omega)\cap H^1(\Omega)$ vanishes on $\partial\Omega$, so
$\partial_iu\in H^1_0(\Omega)$ solves $(\Delta+\mu)\partial_iu=0$, and $\mu<\lambda_D$ gives
$\nabla u=0$, hence $u=0$.

\emph{Case 1: $S_g$ is dense in $\Sone$.} If $\overline{\mathcal B}(S_g)$ were contained in an affine line, so
would be $\overline{\mathcal B}(\Sone)$ by continuity, and a coordinate on that line would contradict
Lemma~\ref{lem:top} with $F=\Sone\setminus A$. So there are three points of $S_g$ with affinely
independent images, hence two indices $j,k$ such that $1,b_j,b_k$ are linearly independent on these
three points: the two independent difference vectors have a nonzero
two-coordinate minor in an orthonormal basis. As $g\ne0$ there, the functions $g,b_jg,b_kg$ are linearly independent, and by
\eqref{eq:traces} so are $u,u_j,u_k$. This contradicts Proposition~\ref{prop:dim}.

\emph{Case 2: $S_g$ is not dense.} Then $g=0$ on an open arc $I_0\subset A$, and $\Phi(I_0)$ is a
regular open boundary arc. Since $\overline{\mathcal B}$ is injective on the nonempty open set $S_g\subset A$, some
$b_j$ is not constant on $S_g$; then $g$ and $b_jg$, and hence $u$ and $u_j$, are linearly
independent. Their gradients vanish on $\Phi(I_0)$ by \eqref{eq:traces}, so both are constant there.
A nonzero combination $w\in V$ has $w=0$ and $\nabla w=0$ on $\Phi(I_0)$. Choose a ball meeting that arc but no corners or other boundary pieces, and
extend $w$ by zero to the exterior portion of the ball. Its zero Dirichlet
trace makes the extension $H^1$, and its zero Neumann trace makes the extension
a weak solution of $(\Delta+\mu)w=0$ throughout the ball: integration by parts
produces no interface term. Interior elliptic regularity then makes this
extension real analytic. It vanishes on an open exterior set, so it vanishes
in the ball and, by interior unique continuation, throughout $\Omega$.
This contradicts the choice of a nonzero $w$. No reflection through a curved
arc, and no analyticity of the boundary, is used.

Hence $\nabla u\ne0$ in $\Omega$. Continuity on the compact closure
guarantees a global maximum and minimum. Neither can occur in the interior,
since an interior extremum would be a critical point. Both are therefore
boundary values, and no interior point can attain either, proving the strict
inequalities.
\end{proof}

\section{Concluding remarks}\label{sec:audit}

The reciprocal-kernel input is proved in Proposition~\ref{prop:compactkernel}
for the bounded compactly supported densities used by the truncation argument,
and Proposition~\ref{prop:cnd} passes to the singular conformal potential. Thus
Theorem~\ref{thm:A} does not assume \cite[Theorem~5.1]{OAI}. The method,
including the regularization and the limiting framework of
\cite[Section~5]{OAI}, is adapted from that preprint, and the variational
framework with curvature is due to Rohleder \cite{Roh21}. The contribution of
this paper is the extension to corners, with the corner and trace estimates it
requires, together with a differential proof of the finite preservation step,
which is an alternative to \cite[Lemma~5.4]{OAI}.

The condition $\alpha_j<1$ is what makes every exponent in
\eqref{eq:gradcorner} positive. At a reentrant corner the first singular
Neumann mode can have an unbounded gradient and the $H^2$ input may fail, and
the method gives no conclusion there. Some domains with a reentrant corner,
such as L-shaped polygons, are treated by a different argument in
\cite{Hat25}, which makes positive use of that singularity. Cusps, infinitely many corners, rough
boundaries and general convex domains are not treated. For the last class one
would need to pass to the limit without losing the strict conclusion, which the
present argument does not provide.

\subsection*{Acknowledgments}
Large language models were used as assistants in developing and checking the
arguments and in preparing the text. The authors have verified the content and
take full responsibility for it.


\clearpage
\begin{thebibliography}{99}\footnotesize
\bibitem{Rau75} J. Rauch, \emph{Five problems: an introduction to the qualitative theory of partial differential equations}, in: Partial Differential Equations and Related Topics, Lecture Notes in Math.~446, Springer, Berlin, 1975, 355--369. \doi{10.1007/BFb0070610}.
\bibitem{BB99} R. Ba\~nuelos and K. Burdzy, \emph{On the ``hot spots'' conjecture of J. Rauch}, J. Funct. Anal. \textbf{164} (1999), no.~1, 1--33. \doi{10.1006/jfan.1999.3397}.
\bibitem{BW99} K. Burdzy and W. Werner, \emph{A counterexample to the ``hot spots'' conjecture}, Ann. of Math. (2) \textbf{149} (1999), no.~1, 309--317. \doi{10.2307/121027}.
\bibitem{dDP24} J. de Dios Pont, \emph{Convex sets can have interior hot spots}, preprint, arXiv:2412.06344 (2024).
\bibitem{JN00} D. Jerison and N. Nadirashvili, \emph{The ``hot spots'' conjecture for domains with two axes of symmetry}, J. Amer. Math. Soc. \textbf{13} (2000), no.~4, 741--772. \doi{10.1090/S0894-0347-00-00346-5}.
\bibitem{AB04} R. Atar and K. Burdzy, \emph{On Neumann eigenfunctions in lip domains}, J. Amer. Math. Soc. \textbf{17} (2004), no.~2, 243--265. \doi{10.1090/S0894-0347-04-00453-9}.
\bibitem{JM20} C. Judge and S. Mondal, \emph{Euclidean triangles have no hot spots}, Ann. of Math. (2) \textbf{191} (2020), no.~1, 167--211. \doi{10.4007/annals.2020.191.1.3}.
\bibitem{JM22} C. Judge and S. Mondal, \emph{Erratum: Euclidean triangles have no hot spots}, Ann. of Math. (2) \textbf{195} (2022), no.~1, 337--362. \doi{10.4007/annals.2022.195.1.5}.
\bibitem{Siu15} B. Siudeja, \emph{Hot spots conjecture for a class of acute triangles}, Math. Z. \textbf{280} (2015), 783--806. \doi{10.1007/s00209-015-1448-1}.
\bibitem{CGY26} H. Chen, C. Gui and R. Yao, \emph{Uniqueness of critical points of the second Neumann eigenfunctions on triangles}, Invent. Math. \textbf{244} (2026), no.~1, 299--353. \doi{10.1007/s00222-025-01398-x}.
\bibitem{DGJYYZ26} H. Deng, C. Gui, X. Jiang, X. Yang, R. Yao and J. Zou, \emph{Symmetry and critical points of second Neumann eigenfunctions on isosceles trapezoids and kites}, preprint, arXiv:2604.19003v4 (2026).
\bibitem{Hat25} L. Hatcher, \emph{The hot spots conjecture for some non-convex polygons}, SIAM J. Math. Anal. \textbf{57} (2025), no.~4, 4588--4608. \doi{10.1137/24M170702X}.
\bibitem{Hat25b} L. Hatcher, \emph{A hot spots theorem for the mixed eigenvalue problem with small Dirichlet region}, J. Spectr. Theory \textbf{15} (2025), 1367--1382. \doi{10.4171/JST/569}.
\bibitem{OAI} OpenAI, \emph{Strict hot spots and absence of interior critical points on smooth simply connected planar domains}, preprint, 24 September 2026. Available in the \texttt{preprints} directory of \url{https://github.com/openai/math}.
\bibitem{Bur26} K. Burdzy, \emph{Removing boundary smoothness from the hot-spots conclusion in the simple-eigenvalue case}, preprint, arXiv:2610.12622 (2026).
\bibitem{Ban26} R. Ba\~nuelos, \emph{Hot spots: from smooth to nonsmooth simply connected planar domains}, preprint, arXiv:2610.32920 (2026).
\bibitem{DJY26} H. Deng, X. Jiang and X. Yang, \emph{Near-boundary rigidity and strict maximum modulus for Neumann eigenfunctions on convex planar domains}, preprint, arXiv:2610.18546 (2026).
\bibitem{Roh21} J. Rohleder, \emph{A new approach to the hot spots conjecture}, preprint, arXiv:2106.05224v4 (2023). \doi{10.48550/arXiv.2106.05224}.
\bibitem{Roh26} J. Rohleder, \emph{A note on hot-spots free subregions of convex domains},
Complex Anal. Oper. Theory \textbf{20} (2026), no.~2, Paper No.~57.
\doi{10.1007/s11785-026-01910-y}.
\bibitem{BP06} R. Ba\~nuelos and M. M. H. Pang, \emph{Level sets of Neumann eigenfunctions},
Indiana Univ. Math. J. \textbf{55} (2006), no.~3, 923--939. \doi{10.1512/iumj.2006.55.2808}.
\bibitem{Fil05} N. Filonov, \emph{On an inequality between Dirichlet and Neumann eigenvalues for the Laplace operator}, St. Petersburg Math. J. \textbf{16} (2005), no.~2, 413--416. \doi{10.1090/S1061-0022-05-00857-5}.
\bibitem{Gri85} P. Grisvard, \emph{Elliptic Problems in Nonsmooth Domains}, Monographs and Studies in Mathematics~24, Pitman, Boston, 1985. Reprinted as Classics in Applied Mathematics~69, SIAM, Philadelphia, 2011, ISBN 978-1-61197-202-3.
\bibitem{Cri18} R. Cristoferi, \emph{A second order local minimality criterion for
the triple junction singularity of the Mumford--Shah functional}, ESAIM Control
Optim. Calc. Var. \textbf{24} (2018), no.~1, 401--435. \doi{10.1051/cocv/2017010}.
\bibitem{Pom92} Ch. Pommerenke, \emph{Boundary Behaviour of Conformal Maps}, Grundlehren der mathematischen Wissenschaften~299, Springer, Berlin, 1992. \doi{10.1007/978-3-662-02770-7}.
\bibitem{PP} K.-M. Perfekt and M. Putinar, \emph{Spectral bounds for the
Neumann--Poincar\'e operator on planar domains with corners}, J. Anal. Math.
\textbf{124} (2014), no.~1, 39--57. \doi{10.1007/s11854-014-0026-5}. See Lemma~4.3 and its proof.
\bibitem{Nad88} N. S. Nadirashvili, \emph{Multiple eigenvalues of the Laplace operator}, Math. USSR-Sb. \textbf{61} (1988), no.~1, 225--238. \doi{10.1070/SM1988v061n01ABEH003204}.
\bibitem{Nag42} M. Nagumo, \emph{\"Uber die Lage der Integralkurven gew\"ohnlicher Differentialgleichungen}, Proc. Phys.-Math. Soc. Japan (3) \textbf{24} (1942), 551--559. \doi{10.11429/ppmsj1919.24.0_551}.
\bibitem{Bre70} H. Brezis, \emph{On a characterization of flow-invariant sets}, Comm. Pure Appl. Math. \textbf{23} (1970), no.~2, 261--263. \doi{10.1002/cpa.3160230211}.
\bibitem{Sch38} I. J. Schoenberg, \emph{Metric spaces and positive definite functions}, Trans. Amer. Math. Soc. \textbf{44} (1938), no.~3, 522--536. \doi{10.1090/S0002-9947-1938-1501980-0}.
\end{thebibliography}
\end{document}